\documentclass[a4paper]{siamltex}       
\usepackage[english]{babel}
\usepackage{graphicx}
\usepackage{amsmath}
\usepackage[misc]{ifsym}
\usepackage{amssymb}
\usepackage{url}

\newcommand{\sigmab}{\boldsymbol\sigma}
\newcommand{\mub}{\boldsymbol\mu}
\newcommand{\deltab}{\boldsymbol\delta}
\newcommand{\bm}[1]{\mathbf{#1}}
\DeclareMathOperator{\Imag}{Im}
\DeclareMathOperator{\Real}{Re}

\begin{document}

\title{Inversion of multi-configuration complex EMI data with
minimum gradient support regularization. \\ A case study
}

\author{Gian Piero Deidda\thanks{Department of Civil and Environmental
Engineering and Architecture, University of Cagliari, via Marengo 2, 09123
Cagliari, Italy, \texttt{gpdeidda@unica.it}.}
\and
Patricia D\'iaz de Alba\thanks{Department of Mathematics and Computer Science,
University of Cagliari, viale Merello 92, 09123 Cagliari, Italy,
\texttt{patricia.diazdealba@gmail.com,rodriguez@unica.it}.}
\and
Giuseppe Rodriguez\footnotemark[2]
\and
Giulio Vignoli\thanks{Department of Civil and Environmental Engineering and
Architecture, University of Cagliari, via Marengo 2, 09123 Cagliari, Italy, and
Groundwater and Quaternary Geology Mapping Department, Geological  Survey of
Denmark and Greenland (GEUS), C.F. M\o llers All\'e 8, 8000  Aarhus C, Denmark,
\texttt{gvignoli@unica.it}.}
}

\maketitle

\begin{abstract}
Frequency-domain electromagnetic instruments allow the collection of data
in different configurations, that is, varying the intercoil spacing, the
frequency, and the height above the ground. Their handy size makes these tools very practical for near-surface characterization in many fields of applications,
for example, precision agriculture, pollution assessments, and shallow geological
investigations. To this end, the inversion of either the real (in-phase) or the
imaginary (quadrature) component of the signal has already been studied.
Furthermore, in many situations, a regularization scheme retrieving
smooth solutions is blindly applied, without taking into account the prior available knowledge. The present work discusses an algorithm for the inversion of the
complex signal in its entirety, as well as a regularization method that promotes the sparsity of the
reconstructed electrical conductivity distribution. This regularization
strategy incorporates a minimum gradient support stabilizer into a truncated
generalized singular value decomposition scheme. The results of the
implementation of this sparsity-enhancing regularization at each step of a
damped Gauss--Newton inversion algorithm (based on a nonlinear forward model)
are compared with the solutions obtained via a standard
smooth stabilizer. An approach for estimating the depth of investigation, that is, the maximum depth that can be investigated by a chosen
instrument configuration in a particular experimental setting, is also
discussed.
The effectiveness and limitations of the whole inversion algorithm are
demonstrated on synthetic and real data sets.
\end{abstract}

\section{Introduction} \label{sec:1}

Frequency-domain electromagnetic induction (EMI) methods have been used
extensively for near-surface characterization
\cite{lsr95,pellerin02,paine03,ovl05,md08,yy10,bgadrc18}. 
The typical measuring device is composed of two electric coils (the transmitter
and the receiver) that are separated by a fixed distance and placed at a known height
above the ground.
The two coil axes are generally aligned either vertically or horizontally with respect
to the surface of the soil.
The transmitting coil generates a primary electromagnetic field $(H_\mathrm{P})$, which
induces eddy currents in the ground, generating a secondary
field $(H_\mathrm{S})$.
The amplitude and phase components of both fields are then sensed by the
receiving coil.
The device stores the ratio between the secondary and the primary fields as a
complex number.

Initially, raw EMI
measurements were used directly for fast mapping of the electrical
conductivity at specific depths, with no time spent on the inversion.
Recent devices are endowed with multiple receivers (multicoil) or use
alternating currents at different frequencies as probe signals
(multifrequency).
Because of their availability, and with the development of efficient inversion
algorithms and powerful computers, EMI data are increasingly
used for reliable (pseudo)
three- and four-dimensional quantitative assessment of the spatial and temporal
variability of the electrical conductivity in the subsurface
\cite{bgadrc18, dcadlrvc18}. These data are usually collected with
both ground-based and airborne systems \cite{mascv12}, and they have begun to be used not only to infer the soil conductivity but also to determine its magnetic
permeability \cite{gslt16, dr16, ddr17}.
The ratio between the secondary and the primary electromagnetic fields provides
information about both the amplitude and the phase of the signal.
The real part (in-phase component) is affected mainly by
the magnetic permeability of the soil, while the imaginary part (out-of-phase
or quadrature component) is affected mainly by its electrical conductivity.
The in-phase or the quadrature component of the signal is typically inverted to reconstruct either the electrical conductivity or the magnetic
permeability of the soil \cite{dfr14,ddr17,dr16}.

In general, an EMI survey consists of many soundings; in the case of
airborne acquisition, for example, there can be hundreds of thousands.
These soundings, measured with multiconfiguration devices at each
specific location, are usually inverted separately and are only a posteriori
stitched together in a (pseudo) two-/three-dimensional fashion. This is still a common
practice even though inversion schemes based on two-/three-dimensional forward modeling are
becoming available and practical for use. However, the advantages of
truly two-/three-dimensional inversion with respect to one-dimensional
approaches are still debatable
\cite{vmac08}. Sometimes, in order to enforce lateral continuity
between the one-dimensional inversion results, the one-dimensional approaches have been extended to incorporate spatial constraints connecting the model parameters from adjacent models \cite{vcas08}.

As in many other fields of application, regularization is usually
performed by imposing smooth constraints. However, this approach is not
always consistent with the true nature of the system under investigation,
as sharp interfaces might be present, for example. In these situations,
a stabilizer selecting the smoothest solution among all the possible solutions
the data, can produce a misleading solution, whereas
compatible with
the data can produce a misleading outcome, whereas
a regularizing term promoting blocky solutions would definitely be
more coherent with the expectations about the target. For these reasons,
over the years, several approaches have been implemented to retrieve
model solutions characterized by sharp boundaries. A particularly
promising strategy is based on so-called minimum gradient support (MGS)
stabilizers \cite{zhdanov02}. This type of stabilizer has been applied to
several kinds of data and has been implemented in diverse inversion
frameworks, for example, inversion of travel-time measurements
\cite{zvu06, vdc12}, electrical resistivity
tomography \cite{fdva15}, and spatially constrained reconstruction of
time-domain electromagnetic data \cite{lvgvmcm15,vfcka15,vsmv17}.
A preliminary application to frequency-domain EMI data was performed by
\cite{ddrv18}.
The MGS stabilizer is a function of a focusing parameter, which influences
the sparsity of the final reconstruction.
Assigning a small value to this variable promotes the presence of blocky
features in the solution, while a large value produces smooth results.

In this work, attention is focused on the inversion of complex-valued
fre\-quency-domain EMI data collected with different configurations by extending
a numerical algorithm discussed by \cite{dfr14, dr16,ddr17}. The
new results are compared to those obtained by inverting the
quadrature component of the signal. Additionally, the
implementation of an MGS-like regularization technique is studied, coupled with the
truncated generalized singular value decomposition (TGSVD) within a
Gauss--Newton algorithm. For a better interpretation of the reconstructed
conductivity, a possible strategy for assesing the depth of
investigation (DOI) is also presented and used.

The paper is structured in six sections. Section\,\ref{sec:2} introduces the
nonlinear forward modeling approach.
In Sect.\,\ref{sec:3}, an inversion algorithm based on the damped
Gauss--Newton method coupled with a TGSVD regularization scheme, which can
process the whole complex signal, is described.
A minimum gradient support stabilizer and a procedure to incorporate it into the
above algorithm is discussed in Sect.\,\ref{sec:3.5}.
To better evaluate the performance of the investigated inversion
strategies, an approach for estimating the
DOI is presented in Sect.\,\ref{sec:4}, leading to the numerical
experiments on synthetic and real data sets reported in Sect.\,\ref{sec:5}.
Section \ref{sec:6} concludes the paper, summarizing its content.

\section{The nonlinear forward model} \label{sec:2}

\begin{figure}[ht!]
\centering \includegraphics[scale=0.5]{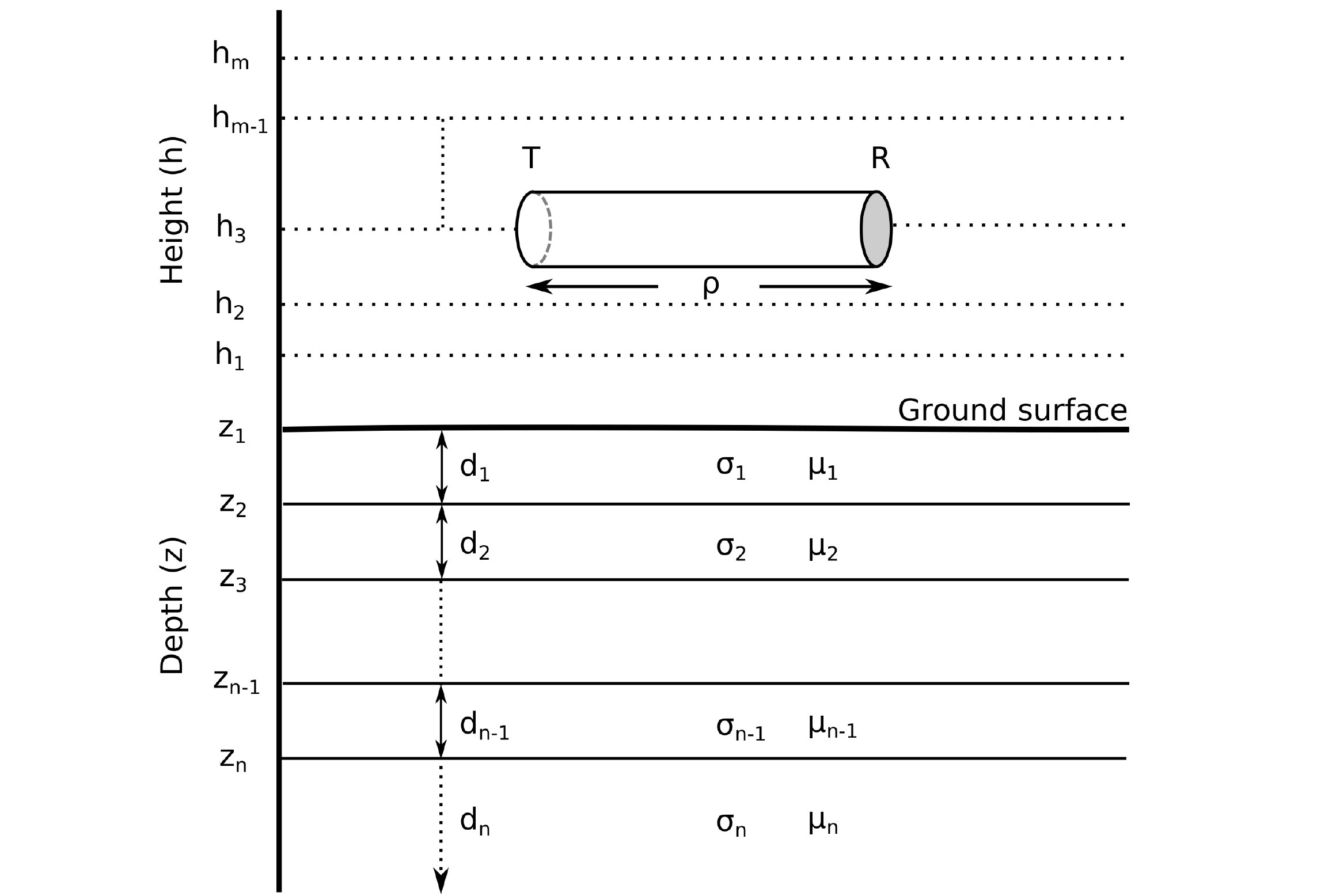}
\caption{Sketch of the subsoil discretization and parameterization along with
the coils of the measuring device above the ground}
\label{fig:0}
\end{figure}

A forward model for predicting the EM response of the
subsoil was discussed by \cite{wait82}. This approach is based on
Maxwell's equations and takes into account the layered symmetry of the
problem. The soil is assumed to have an $n$-layered structure
below ground level ($z_1=0$).
Each horizontal layer, of thickness $d_k$, ranges from depth $z_k$ to
$z_{k+1}$, $k=1,\ldots,n-1$; the deepest layer, starting at $z_n$, is
considered to have infinite thickness $d_n$; see Fig.\,\ref{fig:0}.
The $k$th layer is characterized by an electrical conductivity $\sigma_k$ and a
magnetic permeability $\mu_k$.

The two coils of the EMI measuring device, separated by a
distance $\rho$ and operating at frequency $f$ in Hz, are located at height
$h$ above the ground with their axes oriented either vertically or
horizontally with respect to the ground surface.
Let $\omega = 2\pi f$ be the angular frequency of the device, and let
$\lambda$, ranging from zero to infinity, denote the depth below the ground,
normalized by the intercoil distance $\rho$.
As discussed by \cite{dfr14}, 
if $u_k(\lambda) = \sqrt{\lambda^2+i\sigma_k\mu_k\omega}$ and
$N_k(\lambda) = u_k(\lambda)/(i\mu_k\omega)$ are
the propagation constant and
the characteristic admittance in the $k$th layer, respectively,
then the surface admittance
$Y_k(\lambda)$ at the top of the layer satisfies the recursion equation
\begin{equation}
Y_k(\lambda) = N_k(\lambda)\frac{Y_{k+1}(\lambda)+N_k(\lambda)\tanh(d_k u_k(\lambda))}{N_k(\lambda)+Y_{k+1}(\lambda)\tanh(d_k u_k(\lambda))},
\label{recY}
\end{equation}
for $k=n-1,n-2,\ldots,1$. The recursion relationship in Eq.\,\eqref{recY} is
initiated by setting $Y_n(\lambda) = N_n(\lambda)$ for the deepest layer. It is
worth remarking that both the characteristic and surface admittances depend on
the frequency through the functions $u_k$.

The ratio between the secondary and primary fields for the vertical ($\nu=0$)
and horizontal ($\nu=1$) orientations is given by the expression
\begin{equation}\label{readings}
M_\nu(\sigmab, \mub; h, \omega, \rho) = -\rho^{3-\nu}\int_0^\infty
\lambda^{2-\nu} e^{-2h\lambda} R_{\omega,0}(\lambda) J_\nu(\rho\lambda)\,\mathrm{d}\lambda,
\end{equation}
where $\sigmab = (\sigma_1,\ldots,\sigma_n)^\mathrm{T}$ and $\mub =
(\mu_1,\ldots,\mu_n)^\mathrm{T}$ represent the conductivity and permeability vectors,
respectively, and $J_s(\lambda)$ denotes the first-kind Bessel function of
order $s$ \cite[Sect.\,4.5]{andrews99}.
The reflection factor
\begin{equation}\label{reflfact}
R_{\omega,0}(\lambda) = \frac{N_0(\lambda)-Y_1(\lambda)}{N_0(\lambda)+Y_1(\lambda)},
\end{equation}
can be calculated by setting $N_0(\lambda) = \lambda/(i\mu_0\omega)$ and
computing $Y_1(\lambda)$ via the recursion in Eq.\,\eqref{recY}, where $\mu_0$
represents the magnetic permeability of free space.
The reader should note that the integrand function in Eq.\,\eqref{readings}
depends on the angular frequency $\omega$, as well as on the vectors $\sigmab$
and $\mub$, through the functions $N_0(\lambda)$ and $Y_1(\lambda)$, which
define the reflection factor in Eq.\,\eqref{reflfact}.

The complex-valued functions $M_0$ and $M_1$ can be expressed in a more
compact form in terms of the Hankel transform \cite[Sect.\,4.11]{andrews99}

\begin{equation*}
\mathcal{H}_\nu[f](\rho) = \int_0^\infty f(\lambda) J_\nu(\rho\lambda)\lambda\,\mathrm{d}\lambda,
\end{equation*}
as follows
\begin{equation*}
M_\nu(\sigmab, \mub; h, \omega, \rho) = -\rho^{3-\nu}
\mathcal{H}_\nu[\lambda^{1-\nu}\, e^{-2h\lambda}\, R_{\omega,0}(\lambda)](\rho), 
\qquad \nu=0,1.
\end{equation*}
In general, EMI devices record both the real (in-phase) and imaginary (quadrature) parts of the field ratio.

\section{The inversion scheme} \label{sec:3}

To investigate different depths and be able to infer both the electrical
conductivity and the magnetic permeability profiles for each measurement
location, it is necessary to record EMI data in different configurations. Therefore,
the measurements can be
acquired with different intercoil distances, operating frequencies, and
heights. To further increase the information content in the data,
arbitrary combinations of those configurations can be utilized. Hence, by
indicating with $m_\omega$, $m_h$, and $m_\rho$, respectively, the number of
used frequencies, heights, and intercoil distances, the total number of data
measurements, $b_{tij}^\nu$ (with $t = 1, \ldots, m_\rho$, $i = 1, \ldots,
m_h$, $j = 1, \ldots, m_\omega$, and $\nu = 0,1$) available at each sounding
location is $m =2m_\rho m_h m_\omega$. Of course, the ultimate goal is to
retrieve an estimate of the electrical conductivity vector $\sigmab$
and the magnetic permeability vector $\mub$
that produce the best approximation
$M_\nu(\sigmab,\mub) \approx b_{tij}^\nu$ of the observations.

In the following, it is assumed that the contribution of the
permeability distribution to the overall response is negligible (i.e.,
$\mu_k=\mu_0$ for $k=1,\ldots,n$), so that
the measurements are considered to be sensitive merely to the
conductivity values. However, in principle, the regularization approach
discussed here can also be extended to include the inversion for
the $\mub$ components.
This would require fixing an estimate for the conductivity and determining the
permeability from the data \cite{ddr17},
or computing both quantities by considering the readings defined in
Eq.\,\eqref{readings} as functions of $2n$ variables $\sigma_k$ and $\mu_k$,
for $k=1,\ldots,n$.

To retrieve the conductivity values $\sigma_k$ ($k = 1,\ldots,n$)
associated with the best data approximation, frequency-domain observations
$b_{tij}^\nu$ can be rearranged in a data vector $\bm{b}\in \mathbb{C}^m$; the
same is true for the corresponding calculated responses $M_\nu$, which can be
represented as a vector $\bm{M}(\sigmab)\in \mathbb{C}^m$. Disregarding, for
the moment, the ill-posedness of the problem, the best approximation
$\sigmab^\ast$ can be found by minimizing the Euclidean norm of the residual
$\bm{r}(\sigmab)$, that is,
\begin{equation}
\sigmab^\ast = \arg \min_{\sigmab \in \mathbb{R}^n} \frac{1}{2}\|\bm{r}(\sigmab)\|^2,
\label{min}
\end{equation}
where $\bm{r}(\sigmab) = \bm{b}-\bm{M}(\sigmab)$ takes complex values.

The adopted inversion scheme is based on the Gauss--Newton method, consisting
of the iterative minimization of the norm of a linear approximation of the
residual. Hence, assuming the Fr\'echet differentiability of $\bm{r}(\sigmab)$
\begin{equation*}
\bm{r}(\sigmab_{k+1})\simeq \bm{r}(\sigmab_k)+J_k \bm{q}_k,
\end{equation*}
where $\sigmab_k$ is the current approximation, and $J_k = J(\sigmab_k) \in
\mathbb{C}^{m\times n}$ is the Jacobian of $\bm{r}(\sigmab) =
(r_1(\sigmab),\ldots,r_m(\sigmab))^\mathrm{T}$, defined by $\left[J(\sigmab)\right]_{ij}
= \frac{\partial r_i(\sigmab)}{\partial \sigma_j}$, with $i = 1, \ldots,
m$ and $j = 1, \ldots, n$.

To determine the step length $\bm{q}_k$, as is customary, the real and
the imaginary parts of the arrays involved in the computation are stacked,
\begin{equation*}
\widetilde{\bm{r}}(\sigmab)= \left[\begin{array}{c}
\Real(\bm{r}(\sigmab)) \\
\Imag(\bm{r}(\sigmab))
\end{array}\right] \in \mathbb{R}^{2m}, \qquad
\widetilde{J}(\sigmab)=\left[\begin{array}{c}
\Real(J(\sigmab)) \\
\Imag(J(\sigmab))
\end{array}\right] \in \mathbb{R}^{2m\times n},
\end{equation*}
and the following linear least squares problem is solved:
\begin{equation}
\min_{\bm{q} \in \mathbb{R}^n} \|\widetilde{\bm{r}}(\sigmab_k)
+\widetilde{J}_k \bm{q}\|.
\label{leastsquaresnew}
\end{equation}
The vectors $\widetilde{\bm{M}}(\sigmab),\widetilde{\bm{b}}\in\mathbb{R}^{2m}$
are similarly defined.
In fact, this approach shows that inverting the full complex signal
doubles the number of available data measurements.

The analytical expression of the Jacobian was derived by \cite{dfr14, ddr17}.
In the same papers, it was proven that such an expression is both more
accurate and faster to compute than its finite difference approximation.

To ensure the convergence, while also enforcing the
positivity of the solution, the Gauss--Newton scheme can be
implemented by incorporating a damping factor.
The iterative method becomes
\begin{equation}\label{gniter}
\sigmab_{k+1} = \sigmab_k + \alpha_k \bm{q}_k,
\end{equation}
where the step size $\alpha_k$ is determined according to the Armijo--Goldstein
principle \cite{bjorck96}, with the additional constraint that
the solution must be positive ($\sigmab_{k+1}>0$) at every iteration. This
choice of $\alpha_k$ ensures the convergence of the iterative method,
provided that $\sigmab_k$ is not a critical point, as well as the physical
meaningfulness of the solution.

The inversion of frequency-domain EMI measurements is known to be
ill-posed \cite{zhdanov02}, so that the linearized
problem in Eq.\,\eqref{leastsquaresnew} is severely ill-conditioned for each
value of $k$.
A strategy to tackle the ill-posedness and find a unique and stable solution
consists of including available physical information in the inversion process
via regularization.


A way to incorporate such a priori information in the process is to couple the
original least squares problem, expressed by Eq.\,\eqref{leastsquaresnew}, with
an additional term, leading to the new minimization problem
\begin{equation}
\min_{\bm{q}\in \mathcal{S}} \|L\bm{q}\|^2, \qquad
\mathcal{S}=\{\bm{q}\in\mathbb{R}^n:\bm{q} = \arg
\min\|\widetilde{J}_k\bm{q}+\widetilde{\bm{r}}_k\|\},
\label{regmin}
\end{equation}
where $L$ is a suitable regularization matrix, which defines the $L$-weighted
minimum norm least squares solution \cite{bjorck96}.
The lower the value of $\|L\bm{q}\|$ at the
selected model, the better the matching between the solution and the
a priori information.
By far, the most commonly used regularization matrices favor solutions that are
smoothly varying (either spatially or with respect to a reference model). In
such cases, $L$ is often chosen to be the identity matrix or a discrete
approximation of the first or second spatial derivative.

To cope with the ill-conditioning of the problem, if $L$ is the
identity matrix, the minimum norm solution of Eq.\,\eqref{leastsquaresnew} at
each iteration of the Gauss--Newton method can be computed by the truncated
singular value decomposition (TSVD) of the Jacobian $\widetilde{J}_k$
\cite{hansen98}.
If $L\in\mathbb{R}^{p\times
n}$, with $p \le n$, is different from the identity matrix, then, assuming the
intersection of the null spaces of $\widetilde{J}_k$ and $L$ to be trivial,
Eq.\,\eqref{regmin} can be solved by means of the TGSVD.
SVD (TGSVD).

The following discussion will be limited to the case $2m \ge n \ge p$,
as the situation characterized by $2m<n$ can be treated in a similar manner.
In this case, the generalized singular value decomposition (GSVD)
\begin{equation}
\widetilde{J}_k = U \left[\begin{array}{ll}
\Sigma & \quad 0 \\
0 & \quad I_{n-p}
\end{array} \right] Z^{-1}, \qquad 
L=V\left[ \begin{array}{ll}
M & \quad 0 
\end{array} \right]Z^{-1},
\end{equation}
where $U\in\mathbb{R}^{2m\times n}$ and $V\in\mathbb{R}^{p\times p}$ have
orthonormal columns, $Z\in\mathbb{R}^{n\times n}$ is nonsingular, and
$\Sigma=\diag[\gamma_1,\ldots,\gamma_p]$, $M=\diag[\xi_1,\ldots,\xi_p]$ are
diagonal matrices with nonnegative entries, normalized so that
$\gamma_i^2+\xi_i^2=1$, for $i=1,\ldots,p$.

The TGSVD solution of Eq.\,\eqref{regmin}, with parameter
$\ell=0,1,\ldots,p$, is then defined as
\begin{equation}
\bm{q}_k^{(\ell)} = \sum_{i=p-\ell+1}^p \frac{\bm{u}_i^T
\widetilde{\bm{r}}_k}{\gamma_i}\bm{z}_i + \sum_{i=p+1}^n (\bm{u}_i^T
\widetilde{\bm{r}}_k)\bm{z}_i,
\label{tgsvdsol}
\end{equation}
in which $\bm{u}_i$ and $\bm{z}_i$ are the columns of $U$ and $Z$, respectively.
Removing the first $\ell$ terms in the first summation of
Eq.\,\eqref{tgsvdsol} eliminates the contribution associated with the smallest
$\gamma_i$. This leads to an approximated solution that is more stable, so
$\ell$ acts as a regularization parameter.
For an implementation of the above discussed inversion algorithm, see \cite{ddflr19}.

At each step of the Gauss--Newton iteration, the regularized minimizer
of Eq.\, \eqref{min} is found by solving Eq.\,\eqref{regmin} through the TGSVD
defined in Eq.\,\eqref{tgsvdsol} for a fixed value of the regularization
parameter $\ell$.
Thus, the solution at convergence
$\sigmab^{(\ell)}$ depends on the specific choice of $\ell$. If a
reliable estimation of the noise level in the data is available, the
regularization parameter can be chosen by means of the discrepancy principle,
which requires that the data fitting must match the noise level in the data.
In contrast, other heuristic strategies can be adopted. One of the most
frequently used approaches is the L-curve criterion \cite{hansen98},
based on the reasonable assumption that the most appropriate
choice for the regularization parameter is the one that guarantees the
optimal trade-off between the best data fitting and the most appropriate
stabilization. A comparison of different strategies for estimating the
regularization parameter was presented by \cite{rr13}.
Clearly, the inverse problem can also be approached in a probabilistic framework;
in this case, the solution consists of a posterior probability distribution
that naturally provides an estimation of the uncertainty of the reconstruction
\cite{kaipio2000,vasic2015}. 
The empirical Bayes method presented by \cite{ggh10} supplies a method for
estimating the regularization parameter, in addition to the overall model
uncertainties.

The forward model $\bm{M}(\sigmab)$, described in Sect.\,\ref{sec:2}, is
strongly nonlinear and nonconvex, so its inversion is rather sensitive
to the starting solution $\sigmab_0$ used to initialize the iterative method
defined by Eq.\,\eqref{gniter}.
In our experience, when the noise in the data is normally distributed
and relatively small, as in the numerical experiments on synthetic data of
Sects.\,\ref{sec:5.1} and \ref{sec:5.2}, any reasonable choice of $\sigmab_0$
converges in general to a solution that may not be the best possible but
still maintains physical significance.
In contrast, when the noise type is consistent with real-world
applications (see Sect.\,\ref{sec:5.3}), an accurate choice of $\sigmab_0$
becomes essential for obtaining meaningful results.
In this paper, the simple procedure of repeating the computation with
a few different constant starting models was adopted, selecting the solution
which produced the minimal residual at convergence.
In the future, we plan to investigate the application of global optimization
techniques \cite{horst2013} in order to reduce the importance of a priori
information for choosing the initial solution.
Such global strategies incur a high computational cost, but they are gaining increasing
popularity
\cite{ggh10} because high-performance
parallel computers are now commonly available.

\section{MGS regularization} \label{sec:3.5}

Both the estimation of the regularization parameter and the choice
of the stabilizing term, which incorporates the available a priori information
on the solution, play an essential role in the accuracy of the final result.
When the solution is known (or assumed) to be smooth, a common choice for
$L$ is the discrete approximation of either the first or second spatial
derivative of the conductivity distribution.
Following the same rationale, to maximize the spatial resolution of
the result, whenever the solution is expected to exhibit a blocky structure, a
stabilizer promoting the sparsity of the computed solution and the retrieval of
sharp interfaces should be used instead.

An example of such stabilizers is
the minimum gradient support (MGS) approach \cite{pz99,zhdanov02,vra17}.
It consists of replacing the term $\|L\bm{q}\|^2$ in Eq.\,\eqref{regmin} with
\begin{equation}
S_\tau(\bm{q}) = \sum_{r=1}^p \frac{\left( \frac{(L\bm{q})_r}{\tau q_r}
\right)^2}{\left( \frac{(L\bm{q})_r}{\tau q_r} \right)^2+\epsilon^2},
\label{sharp}
\end{equation}
where $L$ is a regularization matrix, while $\tau$ and $\epsilon$ are free
parameters.
As can be immediately observed, Eq.\,\eqref{sharp} depends only upon the product
$\tau\epsilon$, so in the following, $\epsilon=1$ is fixed, and only
$\tau$ varies.

\cite{vra17} introduced a generalized stabilizing term that reproduces, for
particular values of two parameters, the $L_2$ and $L_1$ norms, the MGS
stabilizer, and others. The authors showed that for small values of $\tau$
Eq.\,\eqref{sharp} approximates an approach proposed by \cite{last1983} that
minimizes the pseudonorm $\|L\bm{q}\|_0$, that is, the number of nonzero
entries in the vector $L\bm{q}$; see also \cite{pz99}, \cite{vdc12}, and
\cite{wei2018}.
Therefore, the nonlinear regularization term $S_\tau(\bm{q})$ favors the
sparsity of the solution and the reconstruction of blocky features.
If $L$ is chosen to be the discretization of the first derivative
$D_1$, the stabilizer introduced in Eq.\,\eqref{sharp} selects the solution
update corresponding to minimal nonvanishing spatial variation.
This is the origin of the ``minimum gradient'' descriptor.
Its clear advantage is that it can mitigate the smearing and blurring effects
of the more standard smooth regularization strategies.

The parameter $\tau$ determines how each term in Eq.\,\eqref{sharp} affects the
overall value.
In particular, as discussed by \cite{vfcka15}, the model updates with
$$
\left(\frac{(L\bm{q})_k}{\tau q_k}\right)^2<1,
$$
are weakly penalized, as the
corresponding term in Eq.\,\eqref{sharp} is small, while updates with derivatives
larger than the threshold defined by $\tau q_k$ may give a contribution close
to 1.
Thus, the MGS stabilizer penalizes the occurrences of variations larger than the
threshold $\tau q_k$, rather than the magnitude of the variations itself.
This, in turn, favors spatially sparse updates.
The threshold defining when an update is to be considered large enough to be
penalized is dynamically chosen, via the parameter $\tau$, as a fraction of the
actual conductivity update $q_k$.
In conclusion, the MGS stabilizer allows for reconstruction of sharp features
while maintaining the smoothing effect of the regularization $L$ for small
variations in the conductivity updates.

In general, applying the nonlinear regularizing term in Eq.\,\eqref{sharp} to a
linear least squares problem requires a larger computational effort
than the standard first/second derivative approach.
In this case, the least squares problem is nonlinear
itself, so Eq.\,\eqref{sharp} can be treated by the main iterative algorithm: it is linearized at each step of the Gauss--Newton method
by evaluating the terms in the
denominator at the previous iteration $\bm{q}_{k-1}$.
At each step, Eq.\,\eqref{leastsquaresnew} is solved by Eq.\,\eqref{regmin},
replacing $\|L\bm{q}\|^2$ by the approximation
\[
S_\tau(\bm{q}) \approx \|D_{\tau,k}L\bm{q}\|^2,
\]
where $D_{\tau,k}$ is the diagonal matrix with elements
\[
(D_{\tau,k})_{i,i} = \frac{1}{\tau (\bm{q}_{k-1})_r}\left[\left(
\frac{(L\bm{q}_{k-1})_r}{\tau (\bm{q}_{k-1})_r}
\right)^2+\epsilon^2\right]^{-\frac{1}{2}}.
\]
In the numerical simulation described in Sect.\,\ref{sec:5}, the regularization
matrix $L$ is always $D_1$.

Every time the forward model is linear, MGS regularization
leads to a convex problem; this was proved, for example, by
\cite{pz99}.
For nonlinear forward problems, such as the
one discussed in the present study, the further nonlinearity introduced by the MGS stabilizer emphasizes the nonconvex nature of the data fitting problem; see the
discussion at the end of Sect.\,\ref{sec:5.3}.
As already noted, the nonconvexity issue could be approached through global
optimization algorithms \cite{horst2013}, but 
approaches that employ available prior information for the starting model
selection are still of some practical interest for their efficiency.

\section{Depth of investigation} \label{sec:4}

The depth of investigation (DOI) usually refers to the depth below which
data collected at the surface are not sensitive to the physical
properties of the subsurface. In short, the DOI provides an
estimation of the maximum depth that can be investigated from the
surface, given a specific device (in a specific configuration)
and the physical properties of the subsoil. Without a DOI
assessment, it is difficult to judge wether the reconstruction at
depth is produced by the data or is merely an effect of the specific
choice of the starting model and/or the inversion strategy.

One way to assess the DOI can be based on the skin depth calculation,
function of the frequency and the medium conductivity \cite{nm89}.
Alternative methods rely on the study of the variability of the solution
as a function of the starting model. For example, \cite{ol99} discussed
the effectiveness of inverting the data with very different initial half-space
conductivities and subsequently comparing the results to determine the depth
threshold up to which the reconstruction was influenced by the data.

Similar to the strategy by \cite{ca12}, the approach proposed here is based on the integrated sensitivity
matrix, as discussed by \cite{zhdanov02}. Hence, in the following, the DOI
is defined as the depth where, for each individual sounding, the
integrated sensitivity values drop below a certain threshold.
With the aim of studying the sensitivity of the data vector
$\widetilde{\bm{b}}=\widetilde{\bm{M}}(\sigmab)$ to a perturbation vector
$\deltab$, the perturbed data
$\widetilde{\bm{b}}_\delta=\widetilde{\bm{M}}(\sigmab+\deltab)$ is taken into
account.
The linearized version of the problem produces the approximation
\begin{equation*}
\widetilde{\bm{b}}_\delta \approx \widetilde{\bm{M}}(\sigmab) +
\widetilde{J}(\sigmab)\deltab,
\end{equation*}
which implies
\begin{equation*}
\delta \widetilde{\bm{b}} = \widetilde{\bm{b}}_\delta-\widetilde{\bm{b}}
\approx \widetilde{J}(\sigmab)\deltab.
\end{equation*}
Then,
\begin{equation*}
\|\delta \widetilde{\bm{b}}\|^2 = \sum_{i=1}^{2m} (\delta \widetilde{b}_i)^2 =
\sum_{i=1}^{2m} \left(\widetilde{J}(\sigmab)\deltab\right)_i^2,
\end{equation*}
where $\delta\widetilde{b}_i$ denotes the $i$th component of
$\delta\widetilde{\bm{b}}$, $i=1,\ldots,2m$.

Now, assuming $\deltab=\varepsilon\bm{e}_r$, where
$\bm{e}_r\in\mathbb{R}^n$ has zero entries except $(\bm{e}_r)_r=1$,
and denoting by $\widetilde{J}_{i,r}$ the $(i,r)$ entry of the Jacobian
$\widetilde{J}(\sigmab)$, the norm of the perturbation takes the form
\begin{equation*}
\|\delta \widetilde{\bm{b}}\|^2 =
\epsilon^2 \sum_{i=1}^{2m} (\widetilde{J}_{i,r})^2.
\end{equation*}
Then, the integrated sensitivity of the data is defined by
\begin{equation*}
\Sigma_r=\frac{\|\delta \widetilde{\bm{b}}\|^2}{\epsilon^2} = \|\widetilde{J}\bm{e}_r\|^2,
\end{equation*}
where $\widetilde{J}\bm{e}_r$ denotes the $r$th column of the Jacobian
matrix.
This measure represents the relative sensitivity of the data vector to a
perturbation in the conductivity of the ground layer at depth $z_r$.

When $\Sigma_r$ decreases significantly with respect to $\Sigma_1$, that is, when
$\Sigma_r<\eta\Sigma_1$ for a fixed tolerance $\eta$, the recovered
conductivity for the $r$th layer is not strictly related to
the data or, thus, to the physical properties of the subsoil.
Then, the depth $z_r$, at which the reduction $\Sigma_r<\eta\Sigma_1$ occurs
is where the DOI is set.
Evidently, there is some degree of arbitrariness in the choice of the
threshold $\eta$ for the decrease in $\Sigma_r$.


\section{Numerical experiments} \label{sec:5}

Numerical experiments were run on a Xeon Gold 6136 computer running the Debian
GNU/Linux operating system, using a MATLAB software package implementing
the algorithms described in this paper \cite{ddflr19}.
The software is available at the web page
\url{http://bugs.unica.it/cana/software/} as the FDEMtools package.

In the numerical tests illustrated in this section, the electrical
conductivity is determined starting from synthetic and experimental data sets
under the assumption that the magnetic permeability can be approximated by that
of empty space.
The results obtained by processing the quadrature component of
the signal will be compared with those derived from the complex signal in its entirety.
Additionally, the MGS sparsity promoting strategy will be compared
with traditional smooth stabilizers.

\subsection{One-dimensional synthetic data} \label{sec:5.1}

A synthetic data set is generated by representing the conductivity as a function
of depth by the following test functions
\begin{itemize}
\item Gaussian profile: $\sigma_1(z)=e^{-(z-1.2)^2}$,
\item Step profile: $\sigma_2(z)=\begin{cases}
0.2, & \quad z < 1, \\
1, & \quad z \in [1,2], \\
0.2, & \quad z > 2.
\end{cases}$
\end{itemize}

Assuming the magnetic permeability to be that of free space ($\mu=\mu_0$)
and the subsoil to be divided into 60 layers ($n=60$) between $z=0$\,m and
$z=3.5$\,m, the forward model described in Sect.\,\ref{sec:2} is applied
in conjunction with a chosen device configuration to reproduce the instrument
readings.
Since the experimental data studied in Sect.\,\ref{sec:5.3} was recorded
by a CMD Explorer ($\rho =1.48, 2.82, 4.49$\,m; $f = 10$\,kHz),
the instrument readings are constructed according to such configuration, assuming the
measurements were acquired at heights $h = 0.9, 1.8$\,m.
This leads to six readings for each coil orientation ($m_h =2$, $m_\rho=3$,
$m_\omega=1$).

To simulate experimental errors, given a vector $\mathbf{w}$ with normally
distributed entries having zero mean and unit variance, the perturbed data
vector $\widetilde{\bm{b}}_\delta$ is determined from the exact data
$\widetilde{\bm{b}}$ by the formula
$$
\widetilde{\bm{b}}_\delta=\widetilde{\bm{b}}+
\frac{\delta \|\widetilde{\bm{b}}\|}{\sqrt{m}} \mathbf{w}.
$$
This implies that
$\|\widetilde{\bm{b}}-\widetilde{\bm{b}}_\delta\|\approx\delta\|\widetilde{\bm{b}}\|$.
In the computed example, $\delta=10^{-3}$.
The equivalent signal-to-noise ratio (in decibels) is
$$
\text{SNR}_\delta = 10\log_{10}
\frac{\|\widetilde{\bm{b}}\|^2}{\|\widetilde{\bm{b}}-\widetilde{\bm{b}}_\delta\|^2}
= 60\text{dB}.
$$
This noise level is unrealistic in real-world applications, in which the
experimental error may be non-Gaussian and highly correlated.
Here, the aim is to test the performance of the inversion algorithm in an ideal
situation.

For all numerical experiments, the regularization parameter $\ell$ (see
Eq.\,\eqref{tgsvdsol}) is chosen by applying the discrepancy principle, as the
noise is Gaussian and its level is exactly known.

\begin{figure}[ht!]
\includegraphics[scale=0.4]{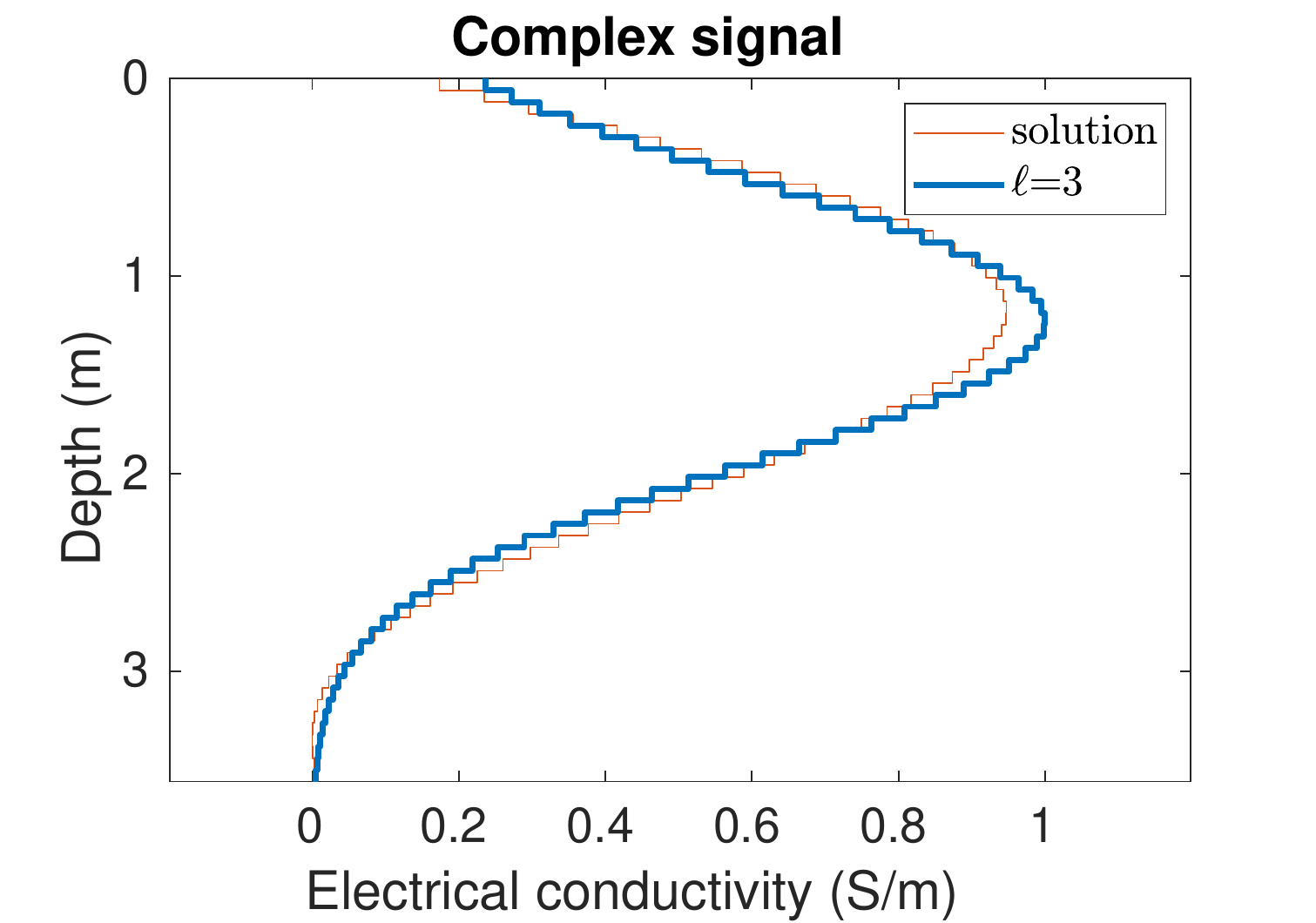} \hfill
\includegraphics[scale=0.4]{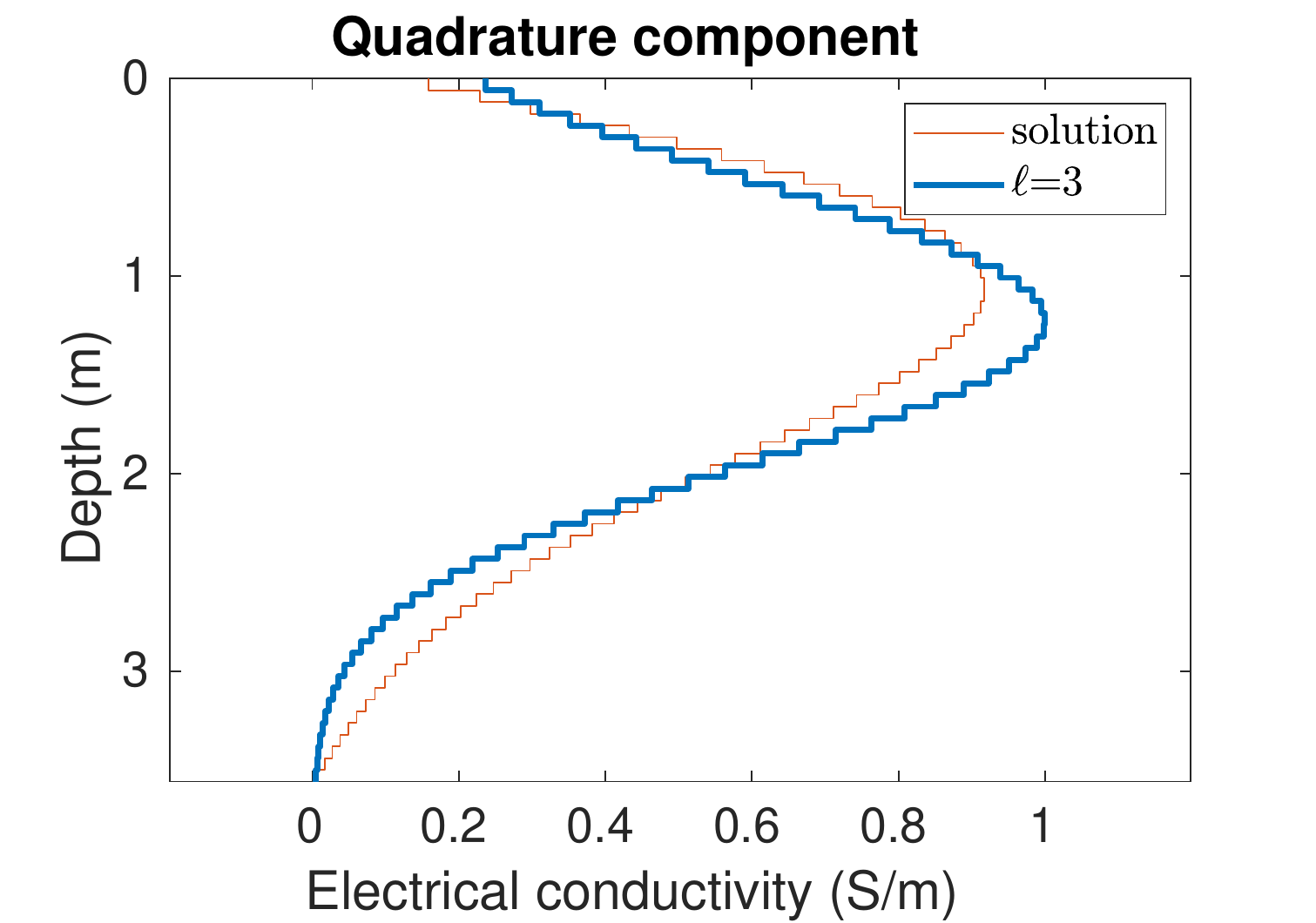} \\
\\
\includegraphics[scale=0.4]{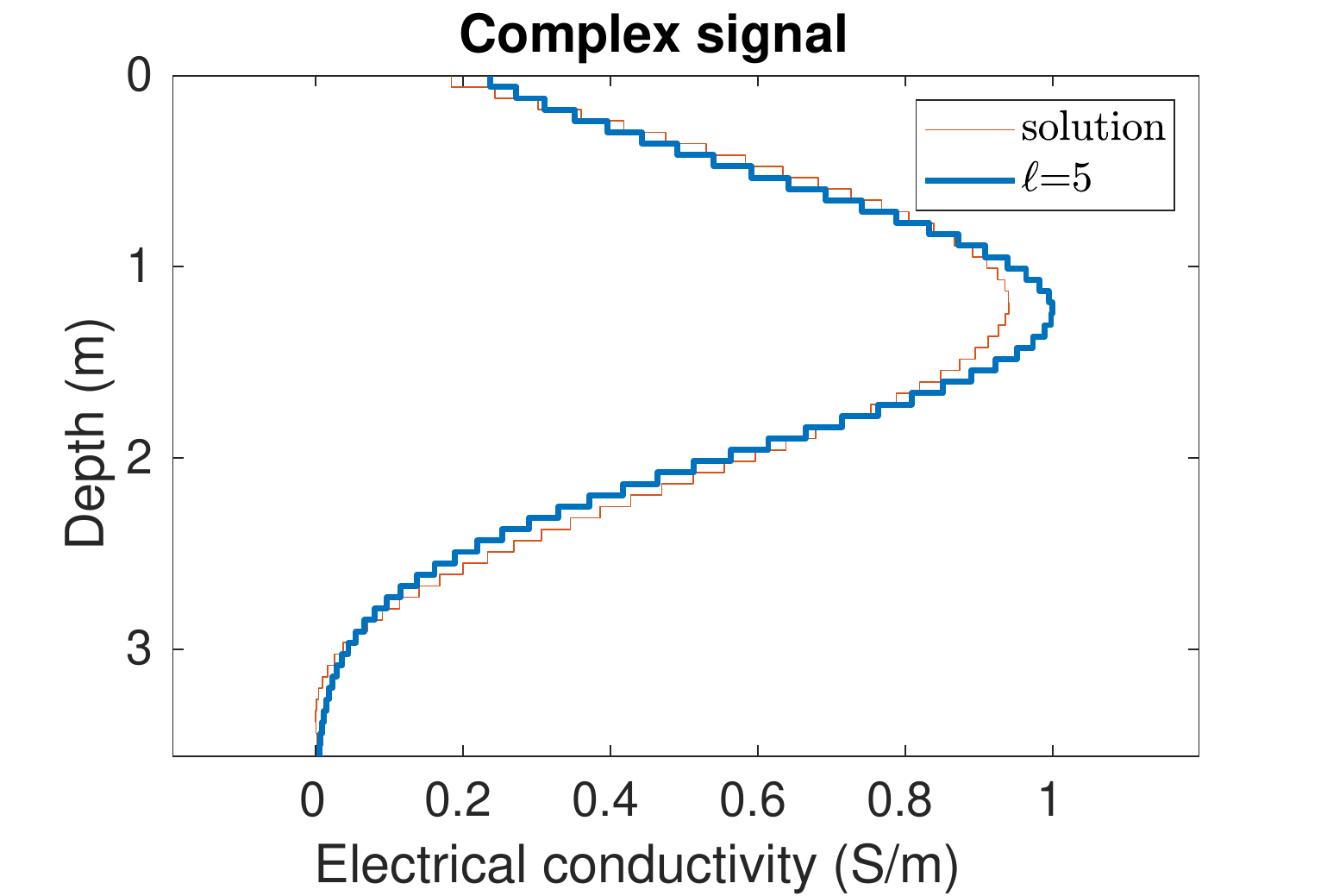} \hfill
\includegraphics[scale=0.4]{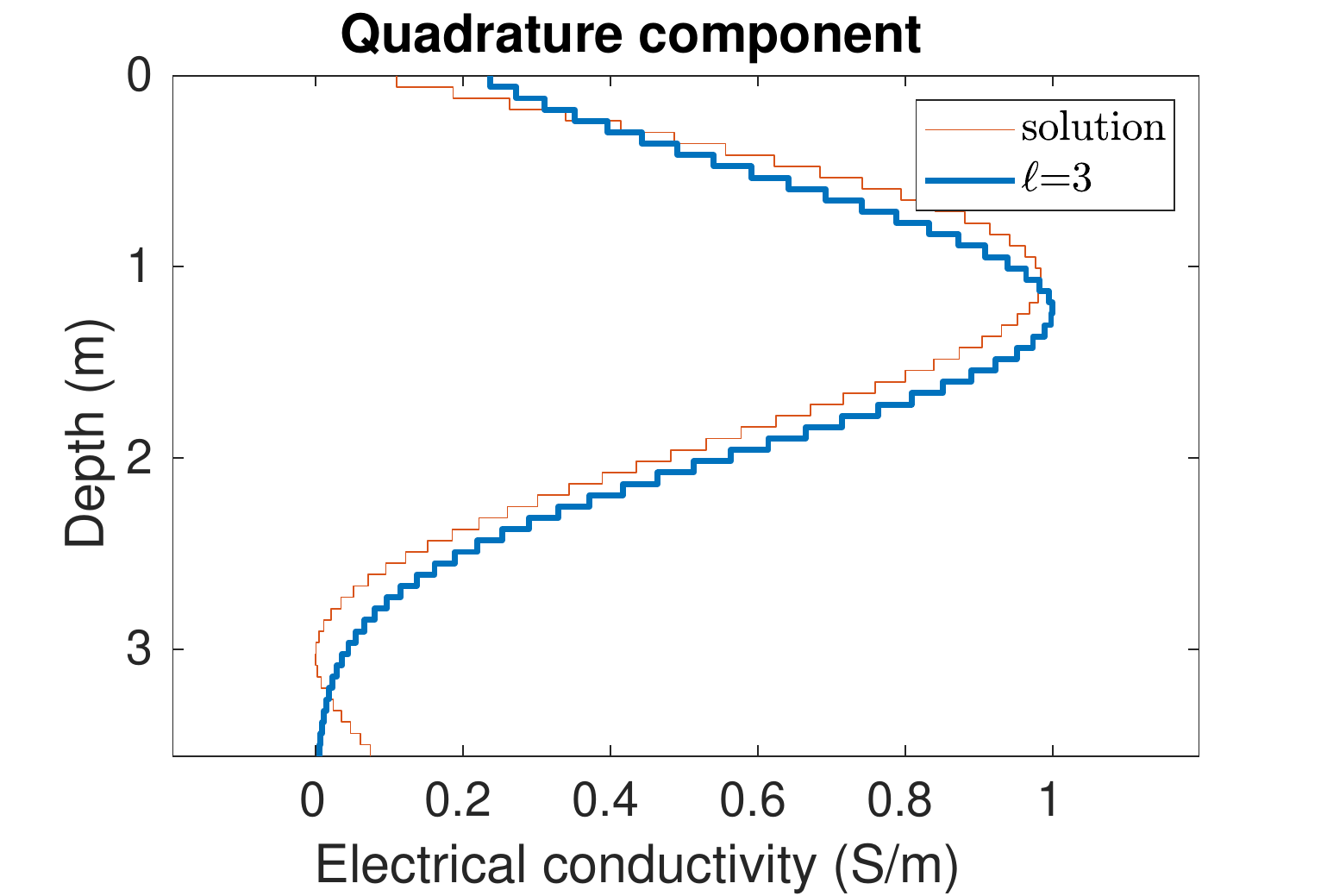}
\caption{Smooth reconstruction of the electrical conductivity for a data set
corresponding to a CMD Explorer configuration with $L=D_2$, test profile
$\sigma_1(z)$, and $\delta=10^{-3}$. Top row: inversion of the complex signal
and of the quadrature component with both coil orientations; bottom row: same
results with the vertical orientation}
\label{fig:1}
\end{figure}

In Fig.\,\ref{fig:1}, the results obtained by the inversion of the
complex signal are compared with those obtained by inverting only the quadrature component.
In this experiment, the smooth test profile $\sigma_1(z)$ and the
regularization term $L=D_2$, the discretization of the second
derivative, are adopted.
The graphs in the top row show the reconstruction of the conductivity when both
orientations of the coils are used, that is, the data set is composed of 12
readings.
The top-left graph represents the solution obtained by inverting the complex
data, while the top-right image reports the reconstruction resulting from
inverting just the quadrature component of the signal.
It is clear that the inversion of the complex signal provides better results.

The graphs in the bottom row of Fig.\,\ref{fig:1} show the results for
the same settings but processing data only for the vertical orientation.
The reconstructions are very similar to those in the top row, showing that
repeating the data acquisition with two orientations of the coils does not
necessarily produce sensibly better results, especially if complex
measurements are processed.

\begin{figure}[ht!]
\includegraphics[scale=0.4]{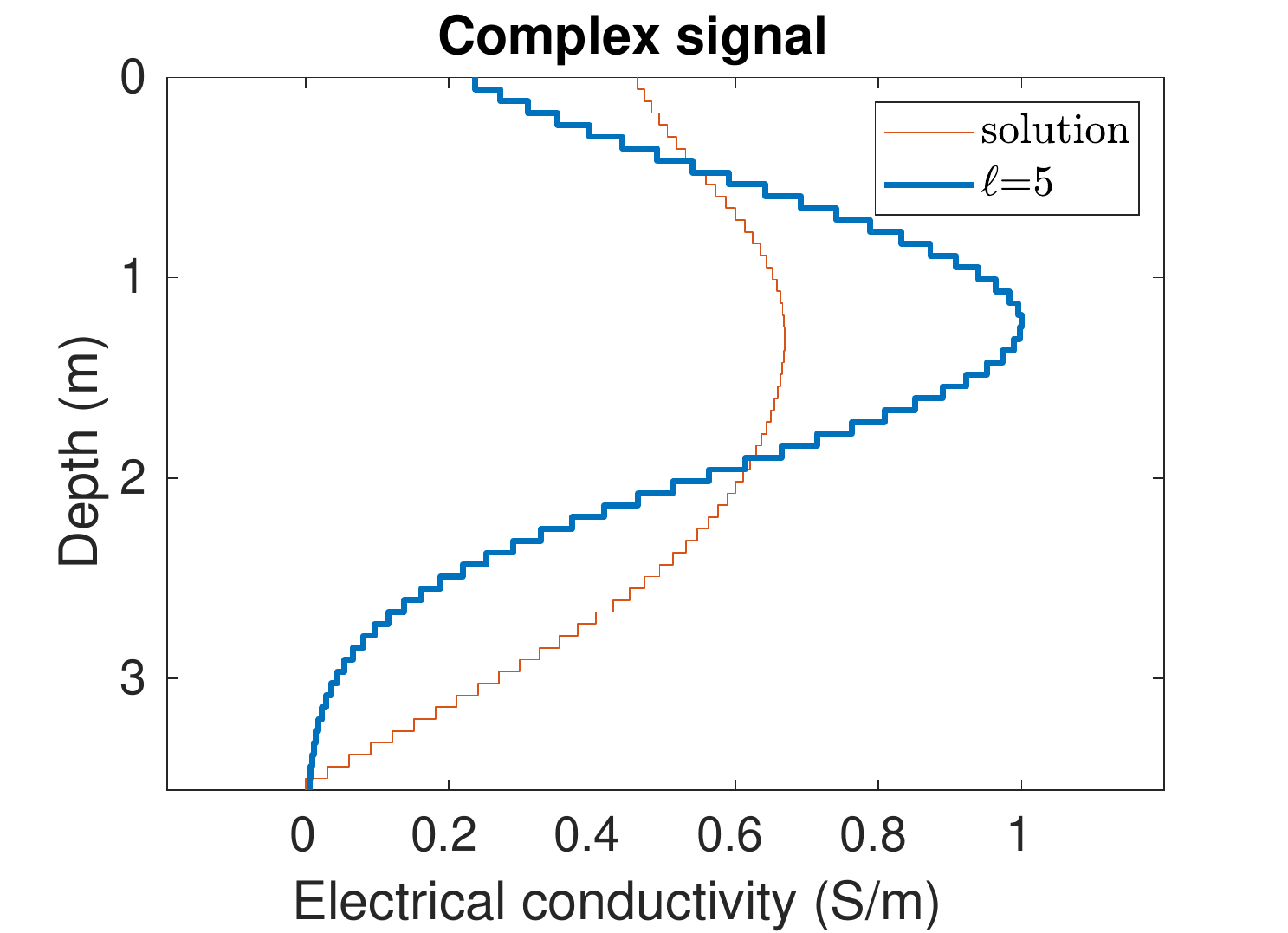} \hfill
\includegraphics[scale=0.4]{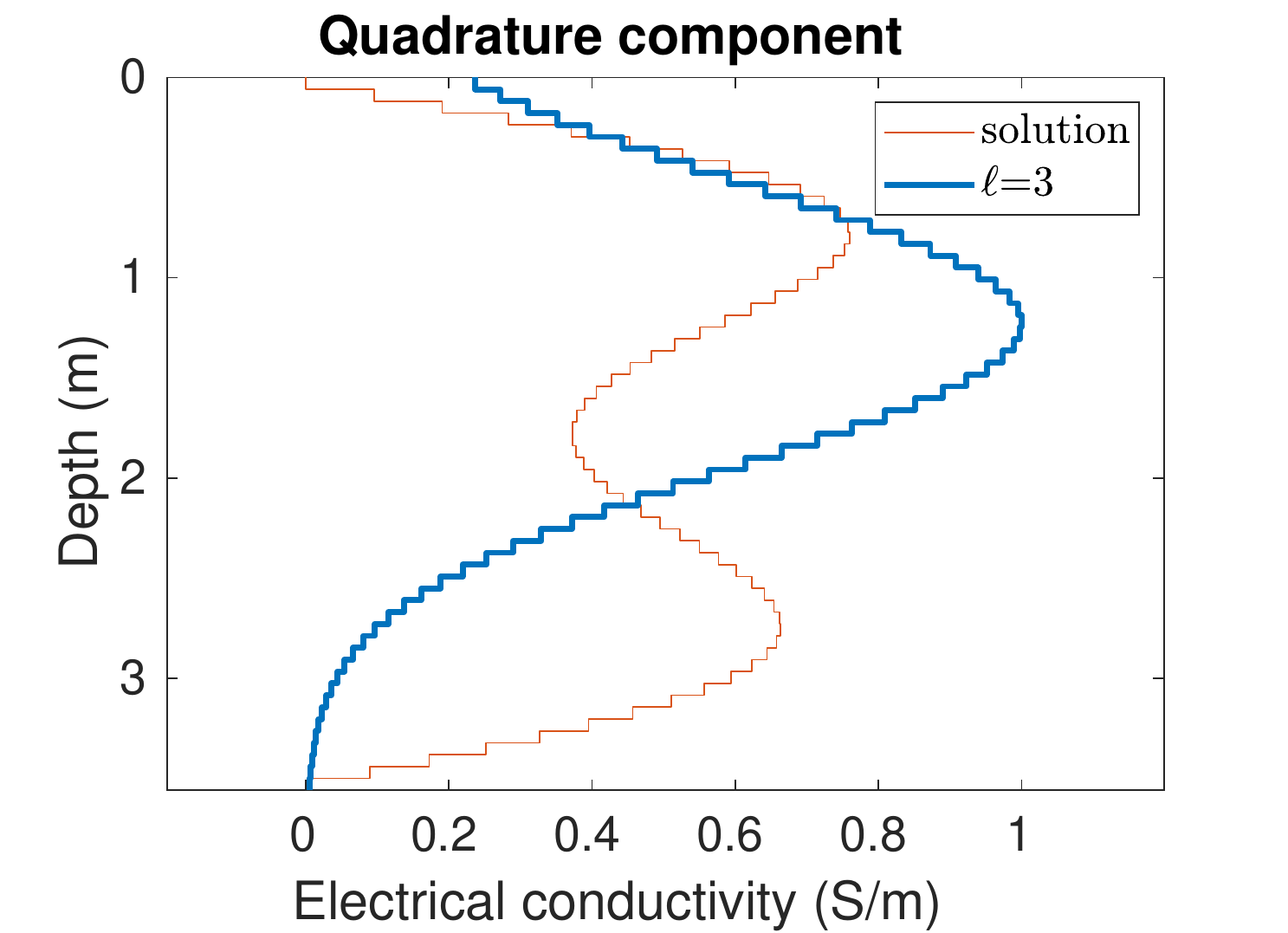}
\caption{Smooth reconstruction of the electrical conductivity for a data set
corresponding to a CMD Explorer configuration with $L=D_2$, test profile
$\sigma_1(z)$, and $\delta=0.2$. The complex signal and of the quadrature
component are inverted with both coil orientations}
\label{fig:1.5}
\end{figure}

To investigate the performance of the algorithm in the presence of
strong noise in the data, the above computation was repeated raising the
Gaussian noise level to $\delta=0.2$, that is, $20\%$ of the signal,
corresponding to $\text{SNR}_\delta=14\,\text{dB}$.
This noise level has been considered realistic in urban sites
\cite{huang2003b}, but it is much larger than that encountered in average
real-world standards \cite{farquharson2003}.
The remaining parameters are kept unchanged, i.e., $L=D_2$,
the discrepancy principle is used to select the regularization parameter
$\ell$, and the 12 readings correspond to both orientations of the coils.
The results, displayed in Fig.\,\ref{fig:1.5}, show that processing the complex
signal leads to reasonably localizing the maximum conductivity.
In contrast, the quadrature component alone does not allow computation of a
meaningful reconstruction.
This example underlines the importance of processing the whole complex
data set when experimental measurements are considered.

\begin{figure}[ht!]
\includegraphics[scale=0.4]{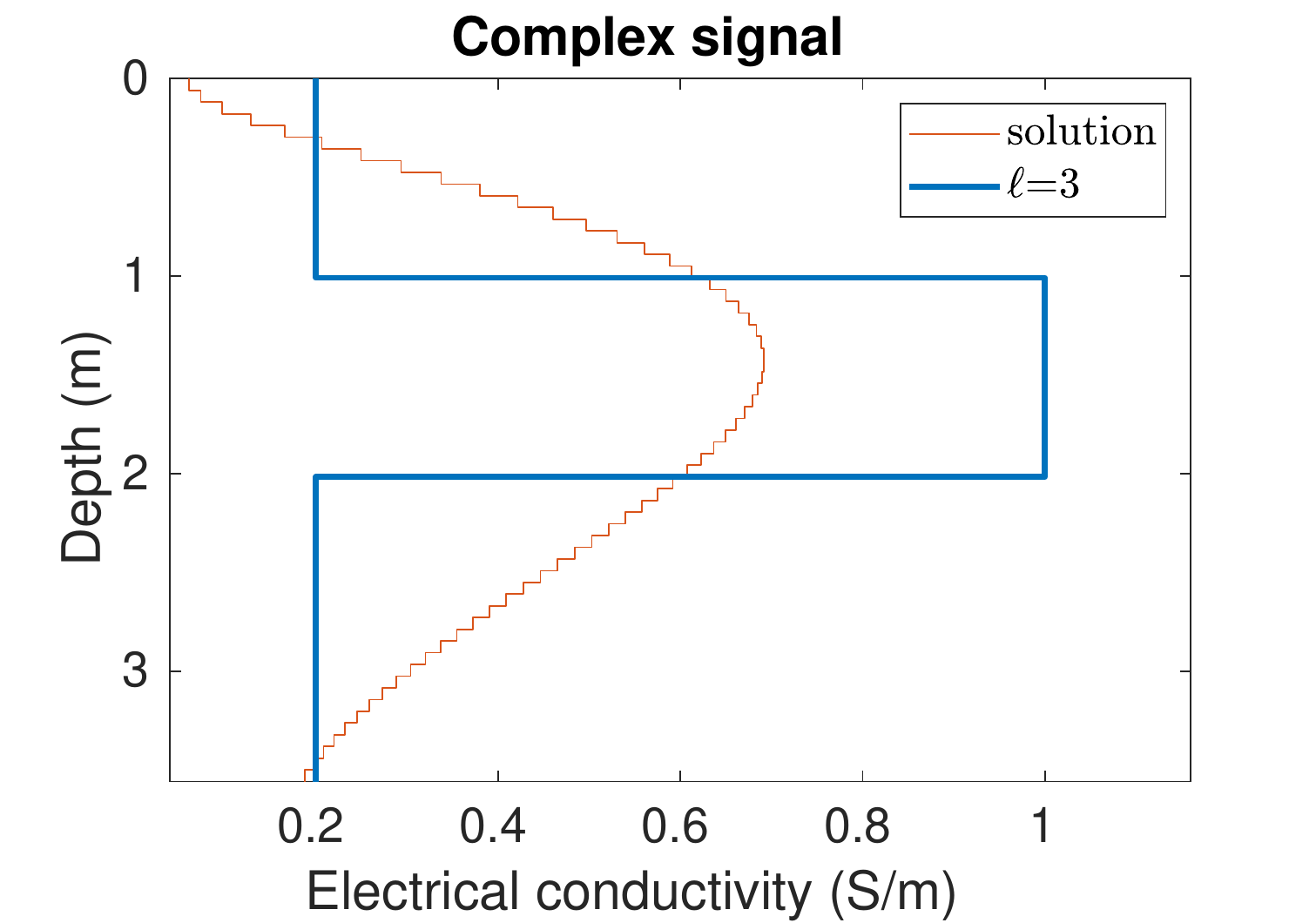} \hfill
\includegraphics[scale=0.4]{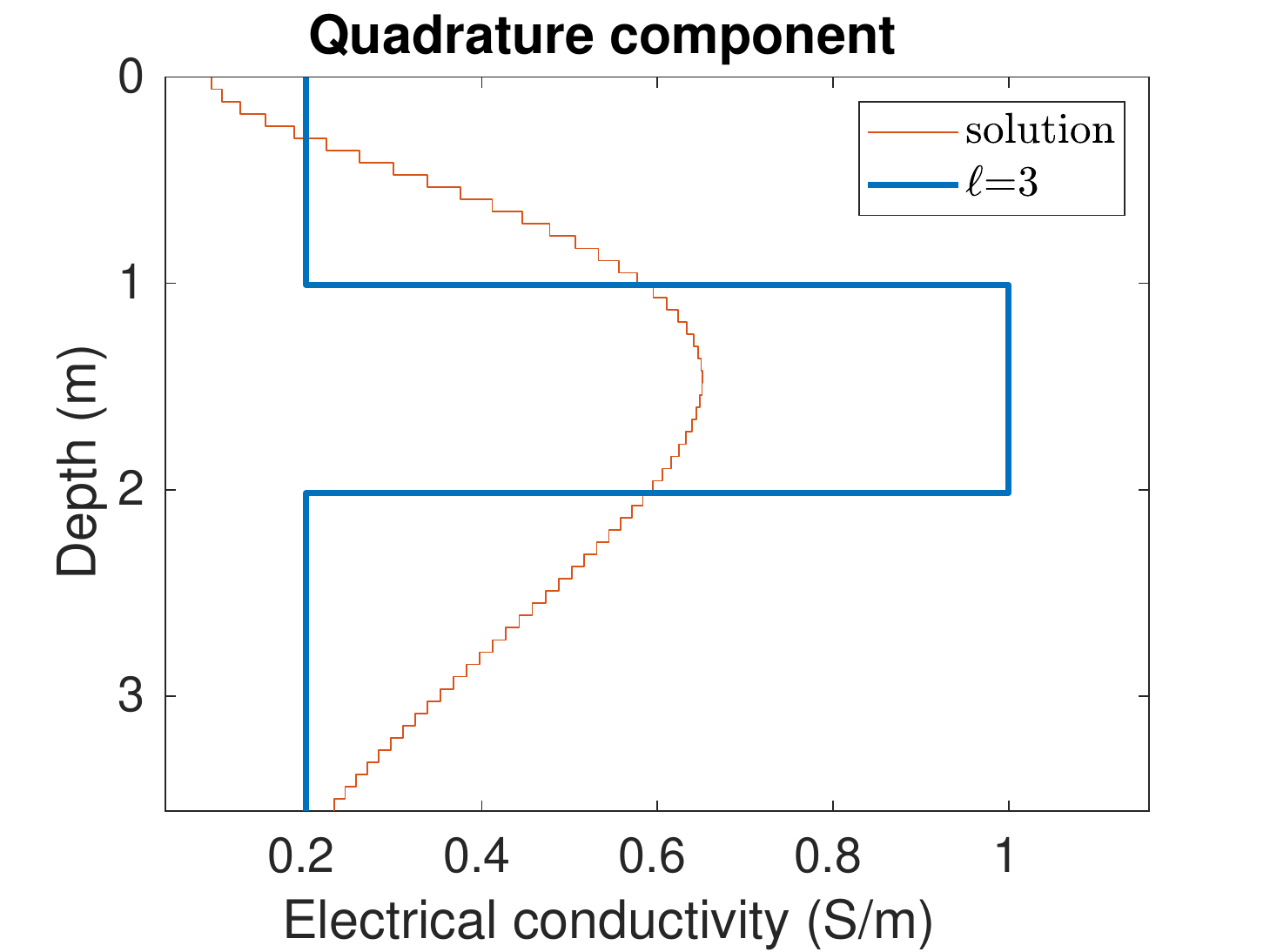} \\
\\
\includegraphics[scale=0.4]{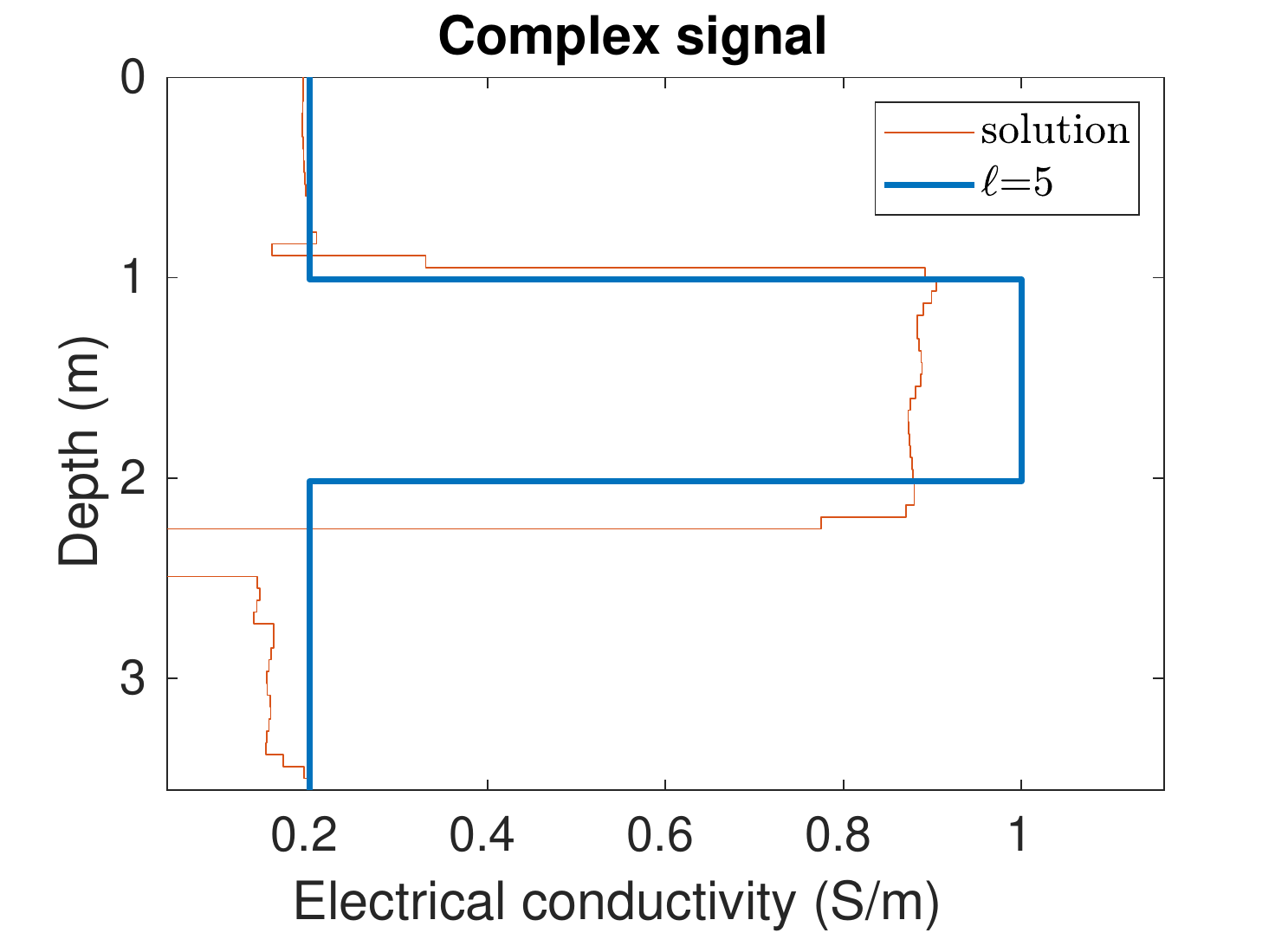} \hfill
\includegraphics[scale=0.4]{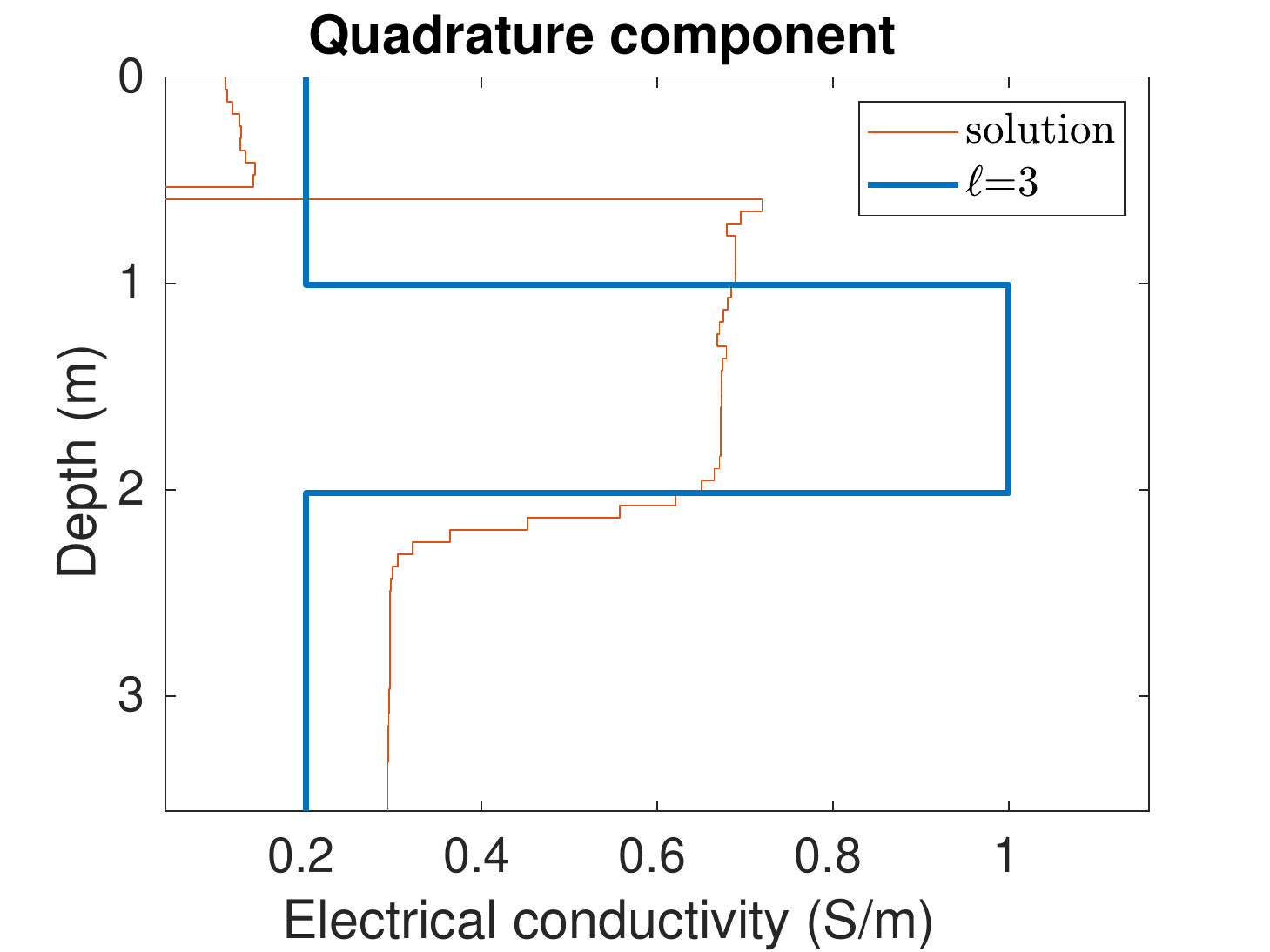}
\caption{Reconstruction of the electrical conductivity for a data set
corresponding to the CMD Explorer configuration and test profile $\sigma_2(z)$,
for $\delta=10^{-3}$ and the vertical orientation of the device. Top row:
smooth inversion of the complex signal and of the quadrature component with
$L=D_1$; bottom row: MGS inversion of the same data}
\label{fig:2}
\end{figure}

Figure \ref{fig:2} displays the results concerning the second synthetic example,
namely, the reconstruction of the discontinuous test profile $\sigma_2(z)$ for
the electrical conductivity.
The same CMD Explorer configuration as before is considered, but the data are
generated only for the vertical orientation of the coils; the noise level is
$\delta = 10^{-3}$.
The graphs in the top row illustrate the performance of the
smooth regularizing matrix $L=D_1$ for both the complex signal and the
quadrature component.
The bottom row displays the same results for the nonlinear regularizing term
$S_\tau(\bm{q})$, after setting $L=D_1$ in Eq.\,\eqref{sharp}.
The results in Fig.\,\ref{fig:2} show that the MGS stabilizer is able to
approximate the presence of sharp boundaries with good accuracy in the model
function. Again, processing the complex signal produces more accurate results.

\subsection{Pseudo two-dimensional synthetic data} \label{sec:5.2}

\begin{figure}[ht!]
\includegraphics[width=\textwidth]{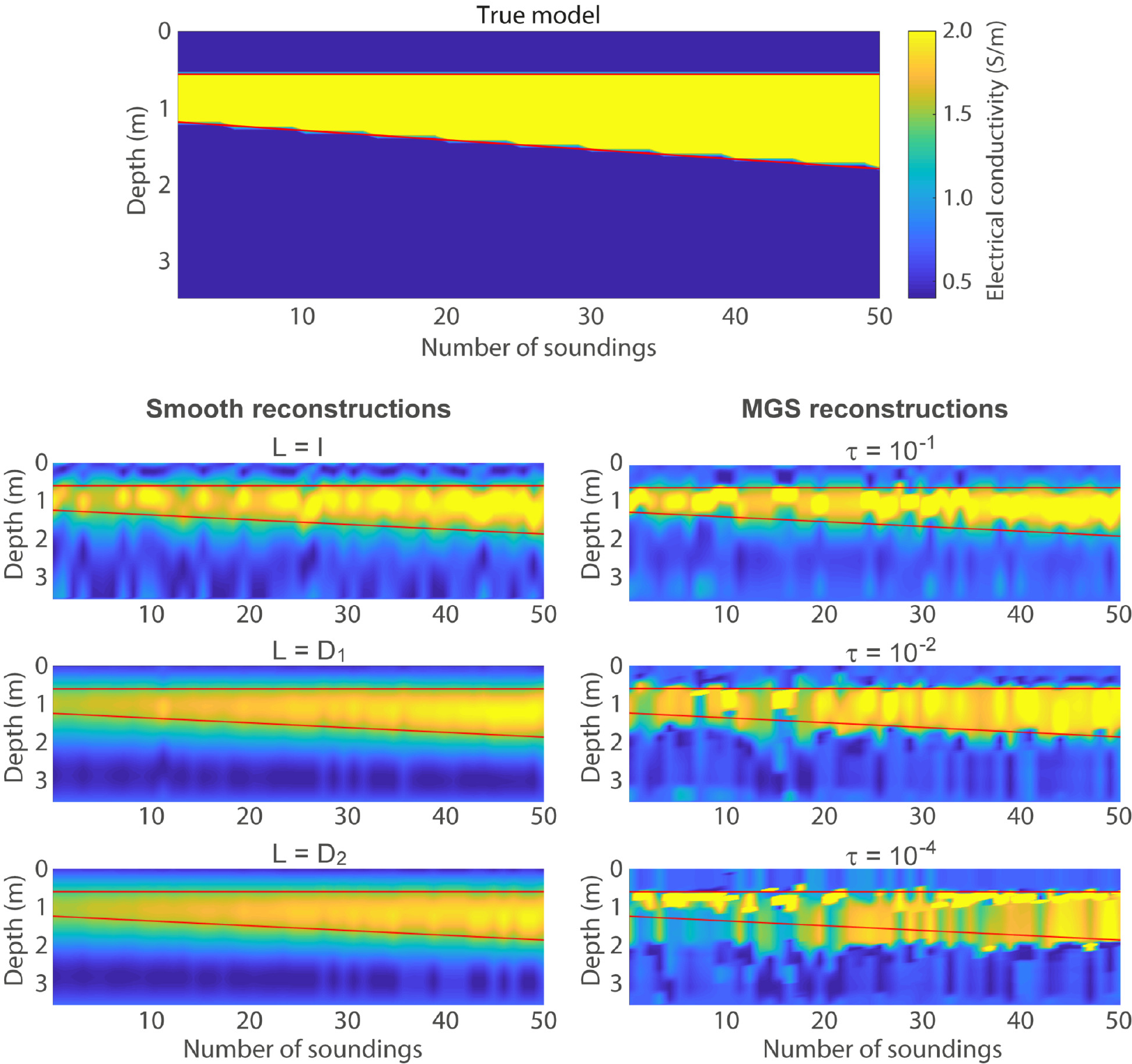}
\caption{Synthetic model for the electrical conductivity (top
graph), smooth inversion of the complex signal (left column), and MGS inversion of
the complex signal (right column)}
\label{fig:3}
\end{figure}

The example described in this section concerns the
reconstruction of a series of one-dimensional models (more
precisely, 50 soundings along a 10\,m straight line) characterized by an abrupt
change in conductivity (from 0.5\,S/m to 2\,S/m) occurring at an increasing depth.
At the top of Fig.\,\ref{fig:3}, the one-dimensional models are depicted side by side in a
pseudo two-dimensional fashion. This facilitates comparison and assessment of the
effectiveness of the methods as the depth of the conductivity transition
varies. The synthetic data simulate an acquisition performed by CMD
Explorer ($\rho=1.48, 2.82, 4.49$\,m, $f=10$\,kHz), with two orientations
of the coils and two measurement heights $h=0.9,1.8$\,m.
The data values are finally perturbed by uncorrelated Gaussian noise with
standard deviation $\delta=10^{-3}$.
To simulate an experimental setting, in which often no information on the noise
level is available, the regularization parameter is estimated in each
one-dimensional inversion by the L-curve criterion.

The left-hand side graphs of Fig.\,\ref{fig:3} show the smooth inversion
results corresponding to $L=I,D_1,D_2$, obtained with a 60-layer
parameterization up to a depth of 3.5 m, with layers of constant thicknesses and a
homogeneous 0.5\,S/m starting model. The right-hand side graphs of
Fig.\,\ref{fig:3} correspond to sharp MGS inversions with three different
values of the focusing parameter $\tau=10^{-1},10^{-2},10^{-4}$.

From these results, it is evident that the smooth inversion for
$L=D_1,D_2$ produces acceptable results, but with an excess of
smoothness. Indeed, it correctly retrieves the transition between the upper
resistivity layer and the lower conductive background, but the transition is
not well identified in space.
It is worth mentioning that the data for each sounding location are generated
independently and that during the inversion, no lateral constraints are
imposed, so the inversion proves to be quite stable.

\begin{figure}[ht!]
\includegraphics[scale=0.4]{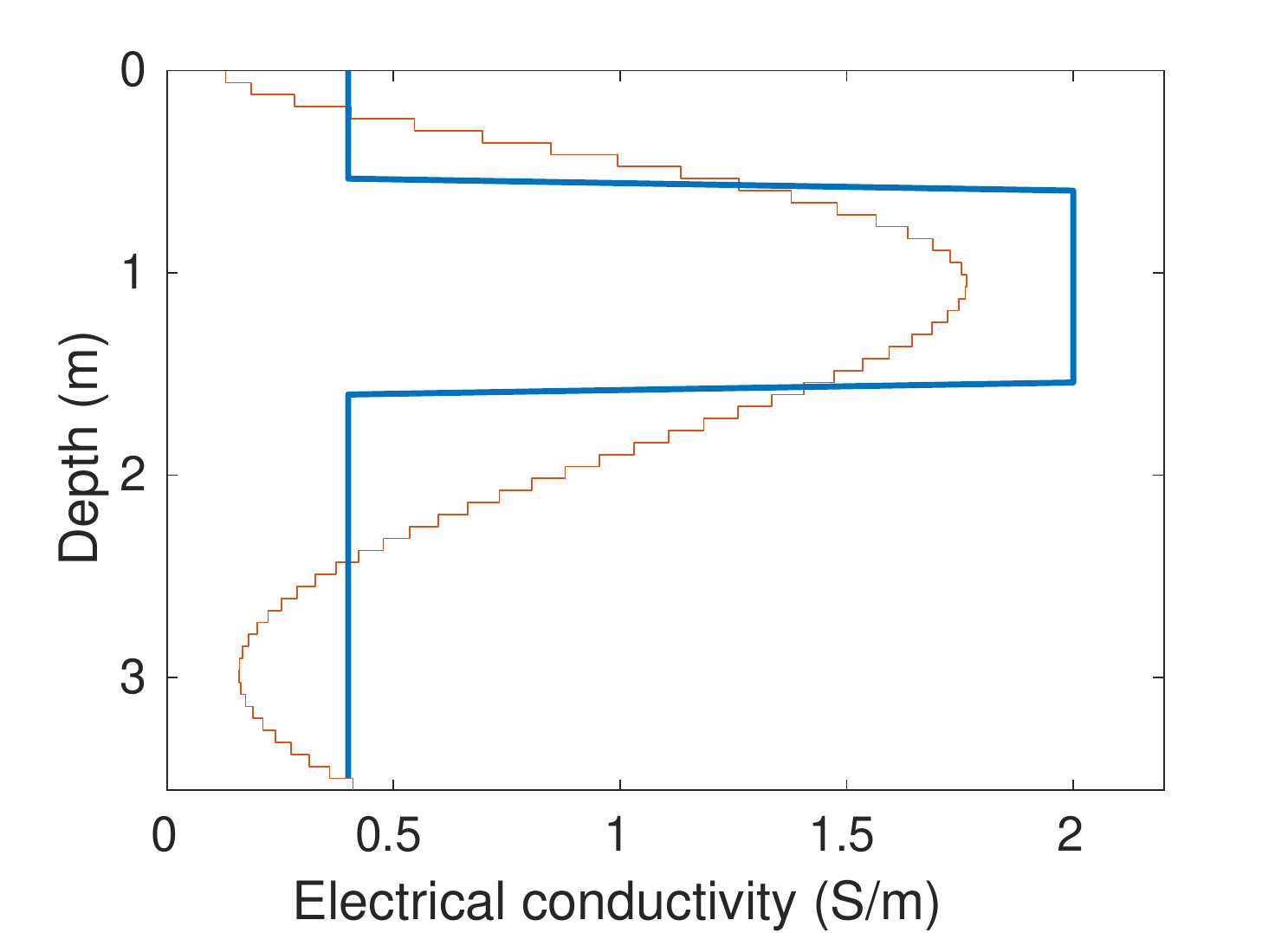} \hfill
\includegraphics[scale=0.4]{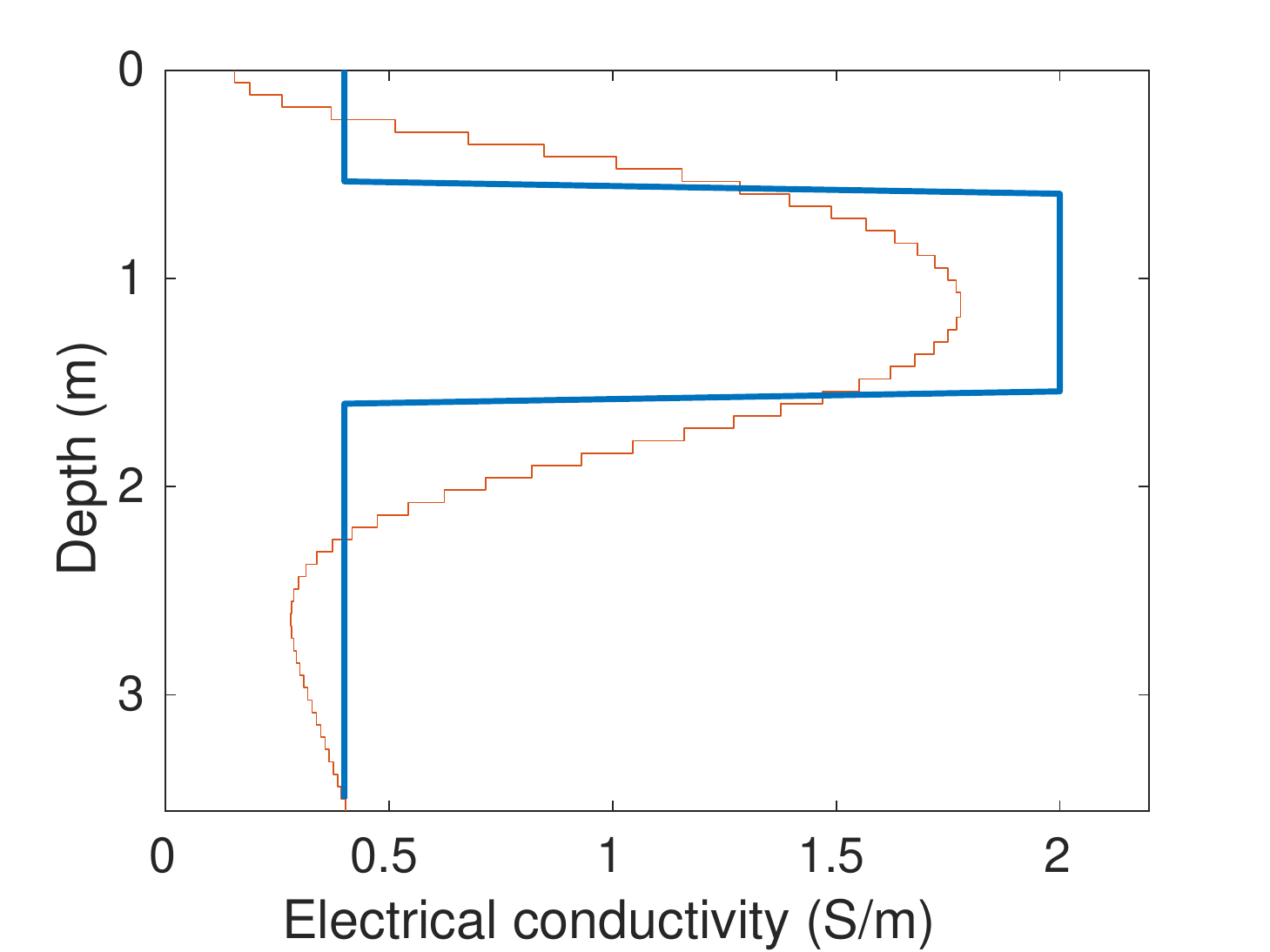} \\
\\
\includegraphics[scale=0.4]{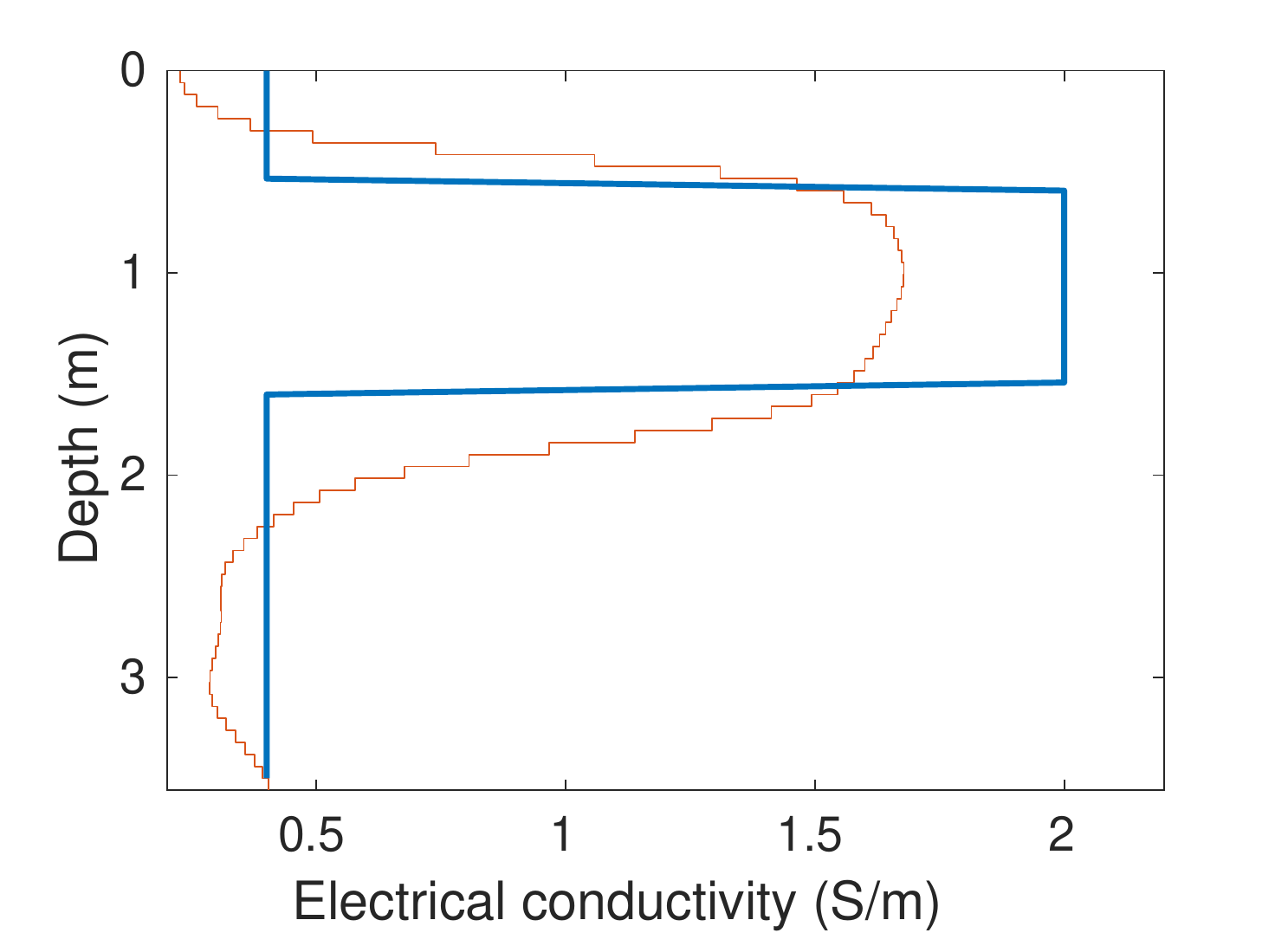} \hfill
\includegraphics[scale=0.4]{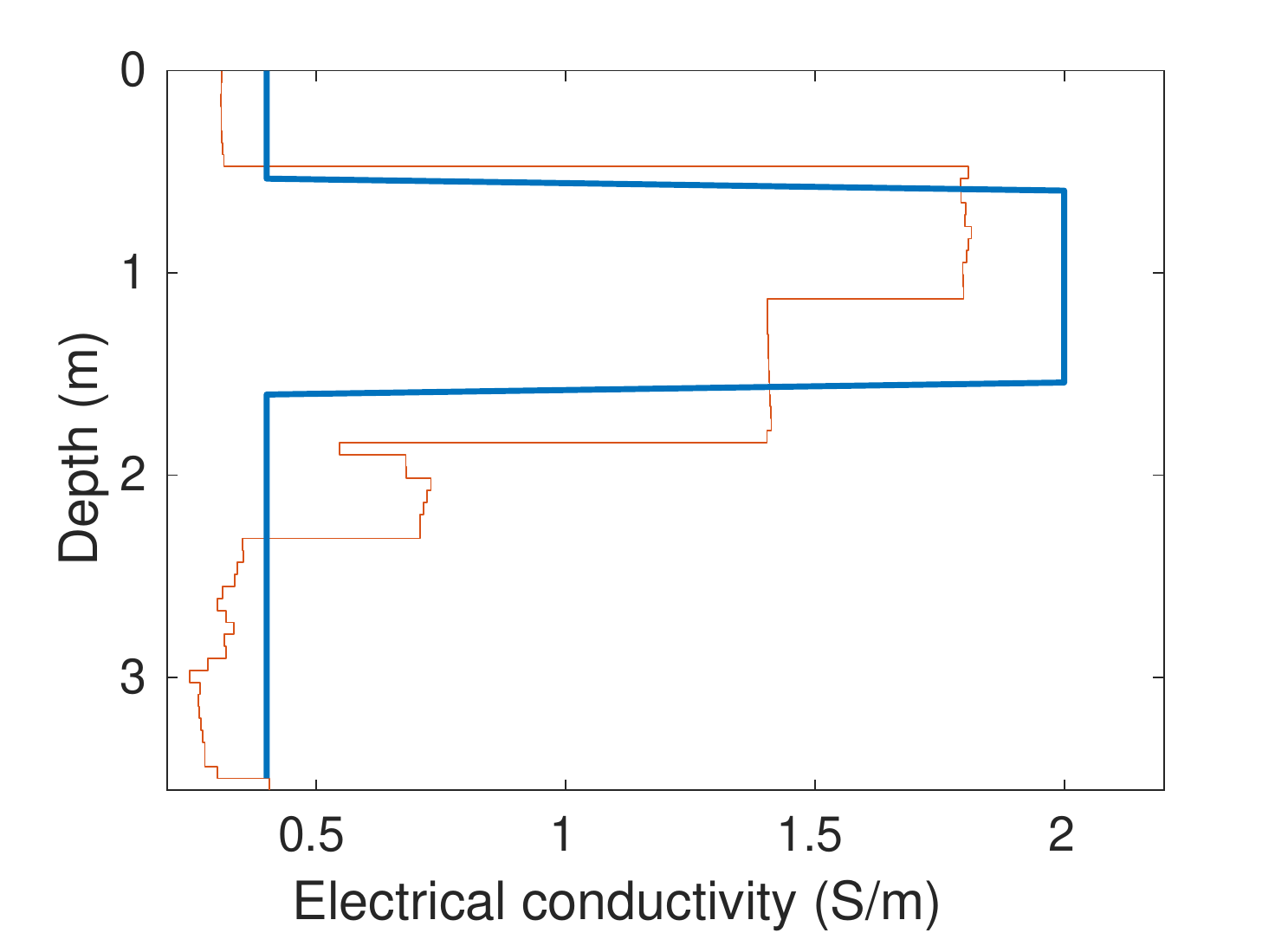}
\caption{One-dimensional reconstruction for the thirtieth column of the
two-dimensional synthetic model of Fig.\,\ref{fig:3}. Top left: smooth inversion
of the complex signal with $L=D_1$; top right: MGS inversion of the
complex signal with $\tau=10^{-1}$; bottom left: MGS inversion with $\tau =
10^{-2}$; bottom right: MGS inversion with $\tau = 10^{-4}$}
\label{fig:4}
\end{figure}

Not surprisingly \cite{vfcka15,fdva15}, the MGS stabilizer with a large $\tau$
produces results very similar to the smooth results.
Decreasing the value of the focusing parameter $\tau$ in Eq.\,\eqref{sharp} corresponds to penalizing the number of small vertical relative variations in the
conductivity updates, as $\tau$ defines the variability range, allowing the derivative update to be considered ``relevant'' for the MGS stabilizer summation.
The sparsity-enhancing effects of the MGS stabilizer are particularly
effective when $\tau=10^{-2}$, where the discontinuity in the solution is more
clearly identified.
When further reducing the focusing parameter, for example, for $\tau = 10^{-4}$,
the reconstructions start to exhibit unrealistic blocky features. This is even
more clear in Fig.\,\ref{fig:4}, where the reconstruction of a single sounding
(the 30th column of the two-dimensional synthetic model of Fig.\,\ref{fig:3}) is
reported, comparing the one-dimensional smooth reconstruction
corresponding to $L=D_1$ (top left) with that of
the MGS stabilizer for $\tau=10^{-1}$
(top right), $\tau=10^{-2}$ (bottom left), and $\tau=10^{-4}$ (bottom right).

A drawback of the MGS inversion process is the instability of the
reconstructed solution. Close one-dimensional reconstructions, that is, close
columns of the ``two-dimensional''
solution, are sometimes very different.
This may be due to the concurrent effects of an incorrect estimation of
the regularization parameter and the nonconvexity of the nonlinear regularized
objective function, leading to solutions becoming trapped in inferior local
optima in some spatial locations.
The reconstruction may be forced to be more regular by imposing
lateral constraints, for example, migrating additional pieces of information
from adjacent models \cite{vfcka15}. This will be the subject of future work.

\subsection{Real Survey} \label{sec:5.3}

\begin{figure}[ht!]
\begin{center}
\includegraphics[scale=0.7]{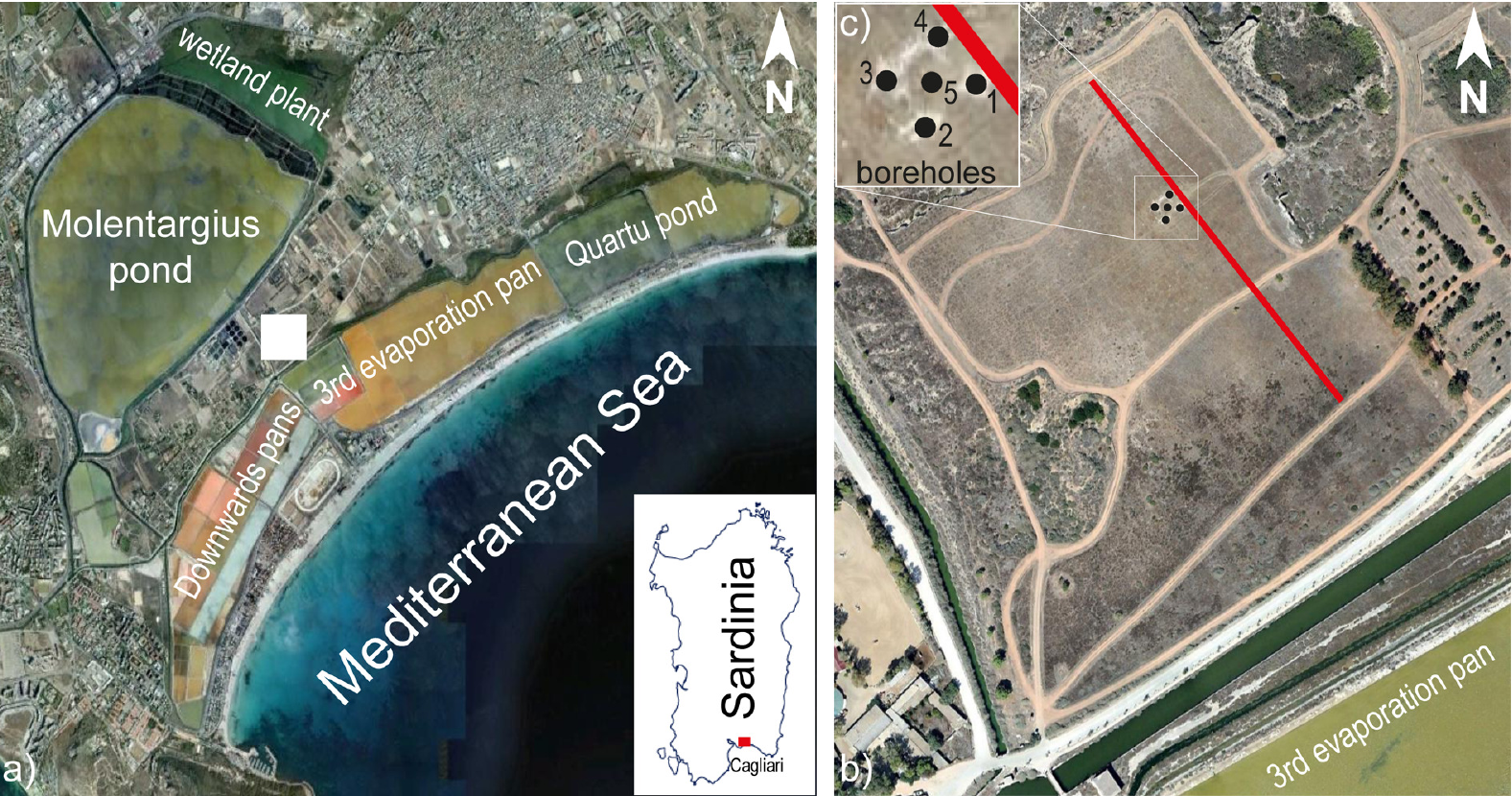}
\caption{\textbf{a} Geographical location of the test site; the white rectangle is the
survey area detailed in (\textbf{b}). \textbf{b} Location of the
electromagnetic profile (red line); black dots indicate the position of the
five boreholes described in \cite{hdcdppsk17}. \textbf{c} Layout and numbering of the
boreholes}
\label{fig:5}
\end{center}
\end{figure}

The proposed algorithm was tested on an experimental data set collected
with a multiconfiguration EMI device at the Molentargius Saline Regional
Nature Park, located east of Cagliari in southern Sardinia, Italy, and
displayed in Fig.\,\ref{fig:5}.
At this site, \cite{hdcdppsk17} investigated the flow
dynamics associated with freshwater injection in a hypersaline aquifer through
hydrogeophysical monitoring and modeling, using five 20\,m deep boreholes
(Fig.\,\ref{fig:5}b, c). The park is a wetland characterized by
the presence of very salty groundwater, with salinity levels as high as 3 times
the NaCl concentration of seawater due to the long-term legacy of infiltration
of hypersaline solutions from nearby salt pans (Fig.\,\ref{fig:5}a) dating
back to Roman times. This site appears to be ideal for testing the MGS
regularized inversion procedure, as the very high electrical conductivity of the
aquifer makes the unsaturated/fully saturated soil interface a sharp electrical
conductivity interface.

\begin{figure}[ht!]
\begin{center}
\includegraphics[scale=0.7]{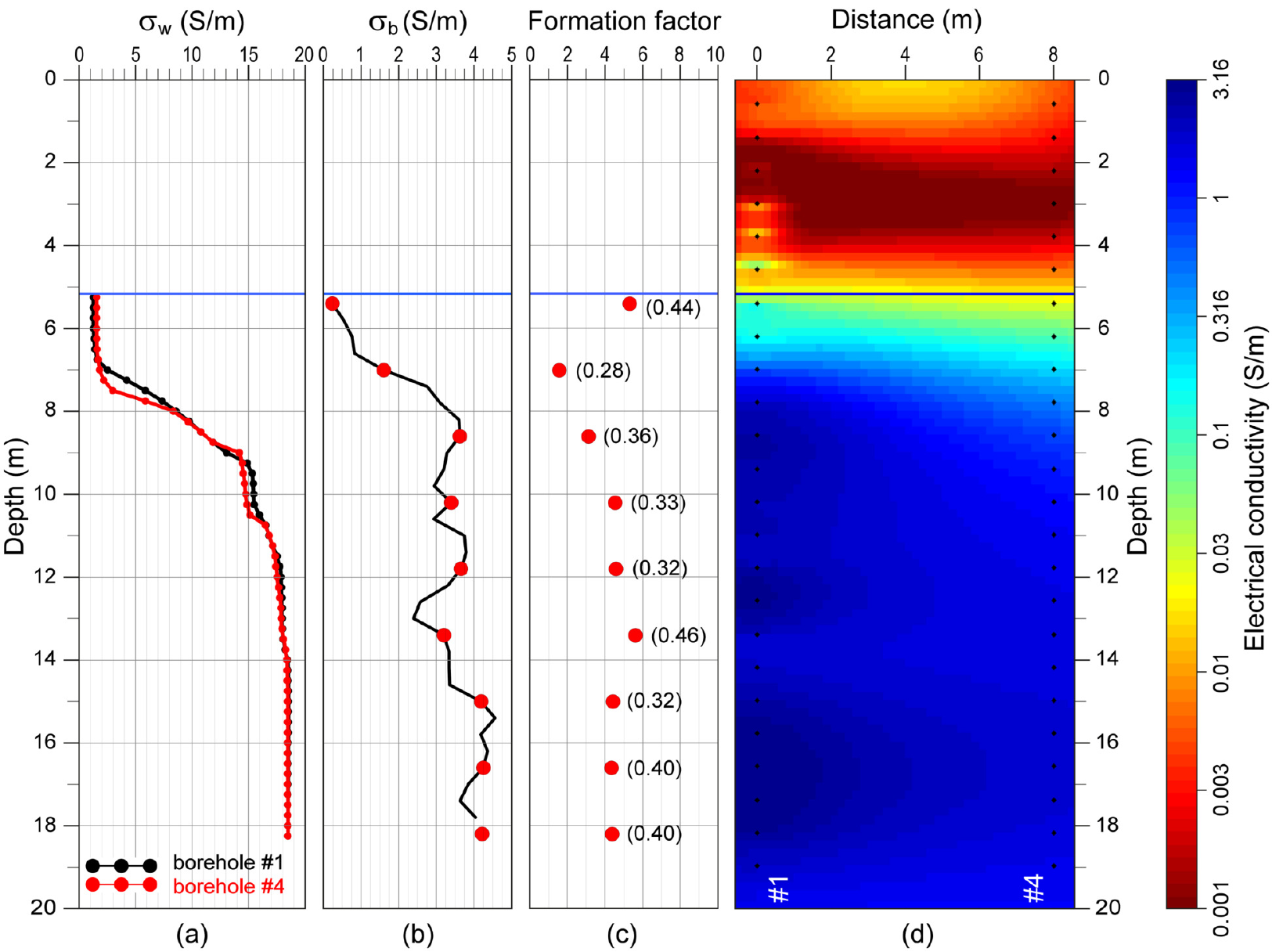}
\caption{\textbf{a} Groundwater conductivity ($\sigma_w$) logs in boreholes \#1 and
\#4; \textbf{b} bulk conductivity ($\sigma_b$) in borehole \#1; \textbf{c} formation factors
along borehole \#1; \textbf{d} cross-hole electrical conductivity image}
\label{fig:6}
\end{center}
\end{figure}

Prior to the freshwater injection experiment, laboratory petrophysical
measurements and different surface, in-hole, and cross-hole electrical
resistivity surveys were carried out to characterize the background of
unsaturated/sa\-tu\-ra\-ted sedimentary succession dominated by sands.
Figure \ref{fig:6} shows some results of these preliminary investigations,
which were used as a reference against which to compare the reliability of the
inversion results.

Groundwater conductivity ($\sigma_w$) logs recorded in boreholes
(see Fig.\,\ref{fig:6}a)
enable two zones to be discriminated, with a transitional 2\,m thick layer in
between: (1) from the water table at 5.2\,m depth to a depth of 7.5\,m, the water
electrical conductivity is approximately 2\,S/m, and (2) below 9.5\,m depth, the water
electrical conductivity reaches 18.5\,S/m.

Using a Terrameter SAS Log (ABEM Instrument)
with an electrode separation of 64
inches, a long normal resistivity log was carried out in borehole \#1 to
estimate the bulk conductivity of the fully saturated soil. Figure
\ref{fig:6}b
shows bulk conductivities $\sigma_b$ calibrated with the values (red dots)
obtained from Archie's empirical relationship \cite{archie1942}
$$
\sigma_b = \frac{\sigma_w}{F},
$$
where $\sigma_w$ is the groundwater conductivity and $F=\phi^{-m}$ is the
formation factor, which is a function of the porosity $\phi$ and the
cementation factor $m$.
The specific values for $F$ were measured from soil samples from borehole \#1
in the laboratory and are shown in Fig.\,\ref{fig:6}c, together with the values
of the porosities of the soil samples (in brackets). These measured bulk
conductivities are apparent conductivities and are representative of a
cylindrical volume with a radius of $\sim 1.5$\,m around the borehole.
They reach values of up to 4\,S/m, but they could be overestimated due to the very high conductivity of the water present in the borehole, which acts as a
preferential path for the current.

Figure \ref{fig:6}d shows the cross-hole conductivity image resulting from the
inversion of apparent electrical conductivities measured with a bipole--bipole
electrode configuration (one current and one potential electrode placed in each
borehole). Black diamonds denote the position of the electrodes, and the blue
line shows the groundwater table at 5.2\,m below the ground surface. Above the
water table, the electrical conductivity is low, ranging between 1 and
10\,mS/m, while in the saturated zone it is very high, and vertical changes due
to the layering of lithologies are not visible. A gradual change to lower
conductivities can be seen only in the upper part, just below the water table.
This is consistent with the water conductivity (Fig.\,\ref{fig:6}a) and bulk
conductivity (Fig.\,\ref{fig:6}b) logs. The conductivity reaches its highest value
below 9.5\,m depth, even if it is slightly smaller than the highest bulk
conductivity, and so it is probably underestimated. This is also an expected
feature for the lack of resolution associated with measurements with large electrode
spacings.

\begin{figure}[ht!]
\begin{center}
\includegraphics[scale=0.7]{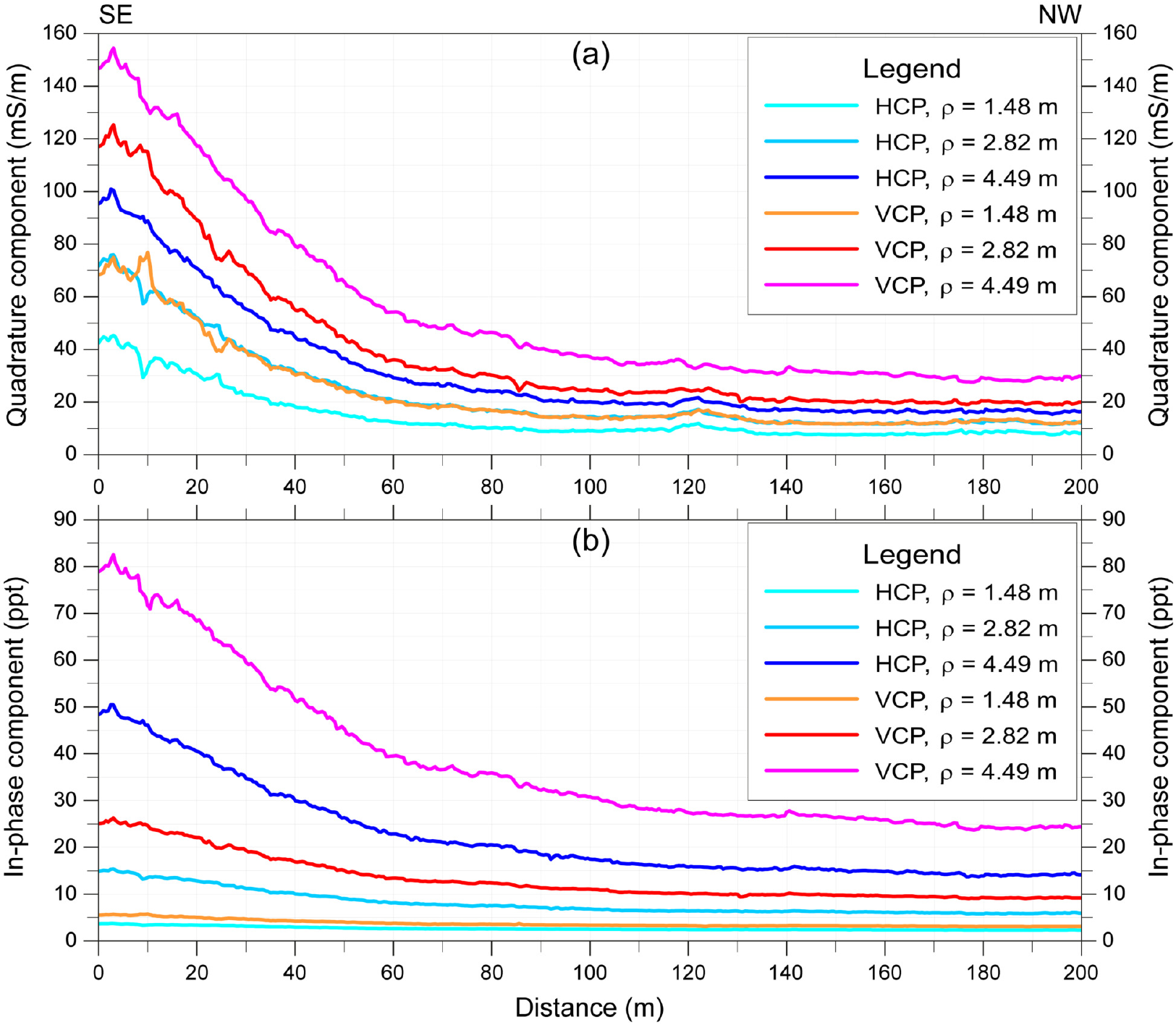}
\caption{EMI raw data recorded along the survey profile: \textbf{a} quadrature
component shown as apparent electrical conductivity; \textbf{b} in-phase component in
parts per thousand (ppt)}
\label{fig:7}
\end{center}
\end{figure}

EMI data were collected along a 200\,m straight-line path (Fig.\,\ref{fig:5}b)
with a topographic elevation varying from 1.6\,m at the southeastern end to
5.7\,m at the northwestern end, using CMD Explorer (Gf-Instruments). This
system operates at a frequency of 10\,kHz and has one transmitter coil paired
with three coplanar receiver coils at 1.48, 2.82, and 4.49\,m from the
transmitter, allowing three simultaneous measurements of the apparent soil
electrical conductivity using vertical (VCP, vertical coplanar) or horizontal
(HCP, horizontal coplanar) dipole configurations. Two surveys were carried out
along the same profile. Data were recorded in continuous mode, with a 0.5\,s time
step, and the system was carried at a height of 0.9\,m above the ground, first
using the HCP and then the VCP dipole configurations. Measurement locations
(UTM coordinates) were logged using a Trimble differential GPS receiver able to
ensure submeter accuracy. Before merging the HCP and VCP data sets, prior to
the inversion, they were spatially resampled at 0.5\,m intervals from a common
starting point to ensure the same number of equally spaced measurement points.
This allowed us to set up a data set consisting of a series of 400 geometric depth
soundings with six complex (quadrature and in-phase components) CMD Explorer
responses each (Fig.\,\ref{fig:7}), suitable for imaging the water table and recovering (by inversion) the
recover (by inversion) the
soil electrical conductivities along the surveyed
profile. At the southeastern end of the survey line, both quadrature and
in-phase responses show higher values than those recorded along the remaining
part, as they were recorded with sensors (transmitter and receiver coils)
closest to the water table.
The data set is available at the web page
\url{http://bugs.unica.it/cana/datasets/}.

The complex response recorded at each sounding point was inverted individually
to infer the electrical conductivity depth profile using the smooth inversion
scheme described in Sect.\,\ref{sec:3}, with the regularization terms $L=D_1$,
$L=D_2$, and the MGS regularization described in Sect.\,\ref{sec:3.5}
with focusing parameter $\tau=10^{-4}$, and $L=D_1$. For all one-dimensional inversions,
the same homogeneous $0.07$\,S/m starting model was used, and the
regularization parameter was estimated by the L-curve criterion.
The discrepancy principle might be used for choosing the
regularization parameter, if a reliable estimation of the noise level were
available. This approach was not pursued in these experiments because our
experience suggests that the noise in EM data is seldom equally distributed
with respect to varying the device configuration; see, for example,
\cite[Fig.\,10]{dfr14}.

Obtaining a noise estimate, even if not essential to perform the computation,
could be useful to better characterize the experimental setting.
Our impression is that the available data set, displayed in Fig.\,\ref{fig:7},
is not sufficient to obtain a trustable noise estimate since repeated
measurements in the same geographical location are missing.
In principle, the linear regularization procedure introduced by
\cite{hrr15} and \cite{prry18} could be adapted to this task.
The method is based on comparing two regularization techniques, for
example, TGSVD and Tikhonov, to select the regularization parameter.
This result is then used to estimate the noise level in the data.
We are currently working on the extension of this method to nonlinear
regularization.

\begin{figure}[ht!]
\begin{center}
\includegraphics[scale=0.7]{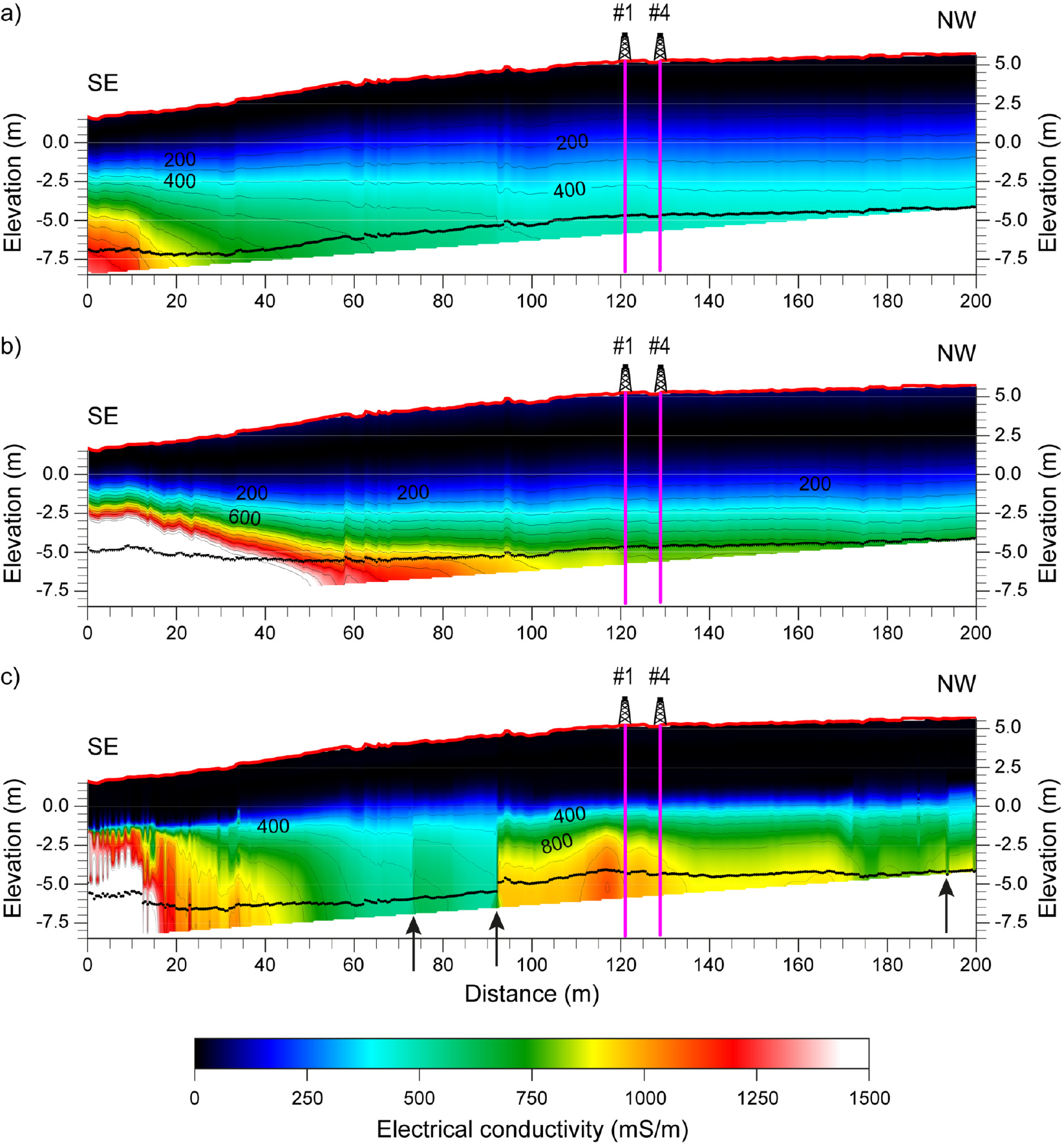}
\caption{Two-dimensional reconstruction of the electrical conductivity from
data collected using CMD Explorer at Molentargius Saline Regional Nature
Park}
\label{fig:8}
\end{center}
\end{figure}

The resulting one-dimensional models,
with 100 layers to a depth of 10\,m below the ground surface
($d_k=0.1\,\text{m}$),
are stitched together and plotted as a pseudo two-dimensional section in Fig.\,\ref{fig:8}.
In each section, the DOI is also plotted to indicate the maximum depth at which the
recovered conductivity is still related to data and not a numerical artifact.
The DOI is represented by the black curves close to the lower boundary of each
section and was estimated by setting the threshold value at $\eta=10^{-2}$ (see
Sect.\,\ref{sec:4}), which was the value that produced results consistent with
the findings of the borehole investigations.

All three solutions satisfactorily capture the overall picture,
although each of them appears better or worse than the
others in certain respects. They clearly retrieve the unsaturated/saturated soil interface at
approximately 0\,m elevation; in the same way, in the southeastern part of the section,
they show the same conductive anomaly due to the saltwater intrusion from the
nearby third evaporation pan of the old saltworks (Fig.\,\ref{fig:5}).

In the smooth solutions (Fig.\,\ref{fig:8}a and~\ref{fig:8}b), the water table
interface is more or less recognizable, but it is not easy to resolve its exact
depth since it does not appear as a sharp interface, as it should actually be
in this case. This undesirable effect, due to the imposed vertical smoothness
constraints, is less noticeable for the solution obtained with the second
derivative regularization term than for the other solutions.
When analyzing the sections in the area between
boreholes \#1 and \#4, indicated by magenta vertical lines, it is clear that
this solution obtains a better fit to the tomography in
Fig.\,\ref{fig:6}d. The solution obtained with $L=D_2$ is also better for the
absolute values of electrical conductivity, which are generally
underestimated when compared with those obtained from the measurements in the
boreholes. Finally, note that in both smoothed solutions, the electrical conductivities vary gradually in the
lateral direction, although they have been obtained by
inverting data, sounding by sounding, without any lateral constraint.

Compared with the previous results, the MGS reconstruction (Fig.\,\ref{fig:8}c) is
less blurred and more reliable in retrieving the sharp water table interface
along the whole section. Electrical conductivities are generally consistent
with those expected on the basis of the results of past surveys. In particular,
to the right of distance 93\,m, the one-dimensional inverted models show
electrical conductivity profiles in very good agreement with those of the
cross-hole tomography. In the MGS solution, however, the electrical conductivities
do not vary gradually in the lateral direction, and they show sharp lateral
changes that do not correspond to real features of the subsoil under
investigation; see, for example, the changes indicated by the arrows in
Fig.\,\ref{fig:8}c.
The reconstruction is particularly erratic in the south-eastern part of the section, where
the nonlinearity of the forward model is amplified by the high conductivity due
to the closeness to the old saltworks; see Fig.\,\ref{fig:5}.

This is again an illustration of the strong dependence of the reconstruction
on the initialization of the iterative method. In these experiments, the same
starting model was adopted for all the data columns. The approach was successful
for the first two regularization matrices, while the MGS stabilizer would
require a more accurate initialization.
This drawback, which was already highlighted at the end of
Sect.\,\ref{sec:3.5}, could be overcome by adopting global optimization
techniques.
Another possibility is to impose a correlation either between the data to be
inverted corresponding to neighboring points or between the obtained
one-dimensional inverse models.

Implementing lateral constraints in the reconstruction, for example, an approach
based on total variation \cite{rudin1992} would couple all the
one-dimensional inversion problems into a large two-dimensional problem,
requiring a suitable solution algorithm.
Indeed, each linearized step of the Gauss--Newton method would lead to a
least squares problem (see Eq.\,\eqref{leastsquaresnew})
too large to be solved by TGSVD.
Another possible approach is based on MGS regularization, followed by
\cite{vfcka15} for another kind of electromagnetic data.

\section{Conclusions} \label{sec:6}

Obtaining relevant solutions to inverse problems requires processing meaningful
data by an effective regularization technique. In the case of EMI data, taking
advantage of both the in-phase and the quadrature component of the available
measurements enriches the data information content, allowing for the
computation of more accurate solutions. Proper regularization consists of the
formalization of a priori information via a stabilizing term. Thus, smoothing
schemes might not always provide the most adequate solution.
Whenever sharp interfaces are expected, it may be wiser to use
regularization terms promoting the sparsity of the retrieved model. In this
manuscript, a Gauss--Newton algorithm regularized by a TGSVD approach,
initially designed for either real (in-phase component) or imaginary
(quadrature component) data inversion, was extended to process complex measurements
and to accommodate an MGS stabilizer. The performance of the new
algorithm was tested on both synthetic and experimental
data sets and compared with alternative approaches.

Synthetic examples over both one-dimensional and pseudo two-dimensional discontinuous
conductivity profiles show that the new algorithm can provide better detail in
the reconstruction than other algorithms.
The MGS solution is not always preferable to the smooth one; it depends on the
expectations/assumptions about the target.
Nevertheless, it is also true that the focusing parameter can be selected such
that the model maintains a certain degree of smoothness.

The enhanced information stemming from complex data values always improved the quality of
the results. The one-dimensional inversion models produced by the complex-valued
experimental data set were able to provide pseudo two-dimensional earth models, consistent
with the findings of in-hole and cross-hole electrical conductivity
investigations, and reliable down to the DOI.

In summary, the new one-dimensional inversion algorithm can be effectively
applied to retrieve smooth and sharp electrical conductivity interfaces in
hydrogeological, soil, and environmental investigations. The instability remains its chief drawback, which may be overcome by adopting global
optimization techniques by developing a more reliable strategy for the
regularization parameter selection and by imposing appropriate lateral
constraints.
This, as well as the application of Bayesian uncertainty quantification ideas,
will be the subject of future work.

\section*{Acknowledgements}
The authors wish to thank the reviewers for their comments, which led to
improvements in the presentation.
This research was supported in part by
the Fondazione di Sardegna 2017 research project ``Algorithms for Approximation
with Applications [Acube]'',
the INdAM-GNCS research project ``Metodi numerici per problemi mal posti'',
the INdAM-GNCS research project ``Discretizzazione di misure, approssimazione
di operatori integrali ed applicazioni'',
the Regione Autonoma della Sardegna research project
``Algorithms and Models for Imaging Science [AMIS]''
(RASSR57257, intervento finanziato con risorse FSC 2014-2020 - Patto per lo
Sviluppo della Regione Sardegna),
the Visiting Scientist Programme 2015/2016 (University of Cagliari),
and an RAS/FBS Grant (Grant No. F71/17000190002).


\end{document}